\documentclass[12pt]{amsart}
\usepackage{amssymb,latexsym,amsmath}
\input xy
\xyoption{all}

\begin{document}

\begin{abstract} Let $\pi$ be a regular algebraic cuspidal automorphic 
representation of ${\rm GL}_n({\mathbb A}_F)$ for a number field $F$. We 
consider certain periods attached to $\pi$. These periods were 
originally defined by Harder when $n=2$, and later by Mahnkopf when $F = 
{\mathbb Q}$. In the first part of the paper we analyze the behaviour of 
these periods upon twisting $\pi$ by algebraic Hecke characters. In the 
latter part of the paper we consider Shimura's periods associated to 
a modular form. If $\varphi_{\chi}$ is the cusp form associated to
a character $\chi$ of a quadratic extension, then we relate the 
periods of $\varphi_{\chi^n}$ to those of $\varphi_{\chi}$, and as
a consequence give another proof of Deligne's conjecture on the critical
values of symmetric power $L$-functions associated to dihedral modular 
forms. Finally, we make some remarks on the symmetric fourth power 
$L$-functions.
\end{abstract}

\title[\bf Periods and special values]
{\bf On certain period relations for cusp forms on ${\rm GL}_n$}

\author{\bf A. Raghuram \ \ \and \ \  \bf Freydoon Shahidi}
\date{\today}
\subjclass{(Primary)11F67. (Secondary)11F70, 22E55}

\maketitle

\tableofcontents

\numberwithin{equation}{subsection}
\newtheorem{thm}[equation]{Theorem}
\newtheorem{cor}[equation]{Corollary}
\newtheorem{lemma}[equation]{Lemma}
\newtheorem{prop}[equation]{Proposition}
\newtheorem{con}[equation]{Conjecture}
\newtheorem{ass}[equation]{Assumption}
\newtheorem{defn}[equation]{Definition}
\newtheorem{rem}[equation]{Remark}
\newtheorem{exer}[equation]{Exercise}
\newtheorem{exam}[equation]{Example}
\newtheorem{quest}[equation]{Question}

\section{Introduction}

The main goal of this paper is to analyze certain periods attached to cuspidal 
representations of ${\rm GL}_n$. If $\pi$ is such a representation, the periods we consider 
are nonzero complex numbers attached to $\pi$ which are (expected to be) intimately linked 
to the special values of the standard $L$-function of $\pi$. In general, these are the only 
periods defined for representations of ${\rm GL}_n$, and ultimately, a study of the special 
values of $L$-functions may be reduced, via functoriality, to a study of these periods.

Let $F$ be a number field and let $\pi$ be a cuspidal automorphic representation of ${\rm 
GL}_n({\mathbb A}_F)$. We assume that $\pi$ is regular and algebraic, which is a condition entirely 
on the infinity component $\pi_{\infty}$ of $\pi$. This assumption makes $\pi$ arithmetically 
interesting, for example, it ensures that the finite part $\pi_f$ is defined over a number field 
${\mathbb Q}(\pi_f)$. Let $S_r$ be the set of real places of $F$. We let $\epsilon = (\epsilon_v)_{v 
\in S_r}$ be a signature indexed by the real places with $\epsilon_v \in \{\pm\}$. 
In the signature, $\epsilon_v$ can be any sign if $n$ is even, and if $n$ is odd, then 
$\epsilon$ is uniquely determined by $\pi$. To this data we 
attach a nonzero complex number $p^{\epsilon}(\pi_f)$ which we call a period of $\pi$. 
See Definition~\ref{defn:period}. 
These periods are defined by comparing a certain canonical ${\mathbb 
Q}(\pi_f)$-structure on the Whittaker model of $\pi_f$ with a ${\mathbb Q}(\pi_f)$-structure on a 
suitable cohomology space to which $\pi$ has nonzero contribution. The comparison map between these 
two spaces is essentially the inverse of the map giving the Fourier coefficients of cusp forms 
in 
the space of $\pi$. These periods were first defined by Harder \cite{harder} for representations of 
${\rm GL}_2({\mathbb A}_F)$, and later were generalized by Mahnkopf \cite{mahnkopf2} in the case of 
${\rm GL}_n({\mathbb A}_{\mathbb Q})$. In both these works they are defined to analyze the special 
values of the standard $L$-function $L(s,\pi_f)$ attached to $\pi_f$.

Concerning special values of $L$-functions, often times, it is interesting to know how these 
values change under functorial operations on the representation at hand. For example, one 
can ask for the behaviour of $L(m,\pi_f \otimes \xi_f)$ as a function of the Hecke 
character $\xi_f$. One application of such a question is the subject of $p$-adic 
$L$-functions. This translates to a question about the behaviour of the periods attached
to $\pi_f$ upon twisting $\pi_f$ by Hecke characters. One of the main aims of this paper
is to answer this question for the periods $p^{\epsilon}(\pi_f)$. In \S\ref{sec:period-twist}
we prove the following 

\smallskip
\noindent{\bf Theorem~\ref{thm:twisted}.}
{\it Let $\pi$ be a regular algebraic cuspidal automorphic representation of 
${\rm GL}_n({\mathbb A}_F)$, 
and let $\xi$ be an algebraic Hecke character of $F$. 
We attach a signature $\epsilon_{\xi}$ to $\xi$. We let $\gamma(\xi_f)$ be the Gauss sum
attached to $\xi$. Then 
$$
p^{\epsilon \cdot \epsilon_{\xi}}(\pi_f) \sim_{{\mathbb Q}(\pi_f,\xi_f)} 
\gamma(\xi_f)^{-n(n-1)/2}p^{\epsilon}(\pi_f)
$$
for any permissible signature $\epsilon$ for $\pi$ (which is an issue only when $n$ is odd).
By $\sim_{{\mathbb Q}(\pi_f,\xi_f)}$ we mean up to an element of the number field ${\mathbb 
Q}(\pi_f,\xi_f)$. Moreover, 
the quotient
$p^{\epsilon \cdot \epsilon_{\xi}}(\pi_f)/(\gamma(\xi_f)^{-n(n-1)/2}
p^{\epsilon}(\pi_f))$ is equivariant under the action of
the automorphism group of complex numbers.}

\smallskip

The proof of this theorem is a little involved to explain in the introduction, however we 
ask the reader to look at the diagram of maps (\ref{eqn:cube}). The proof comes out of an 
analysis of that diagram; the period relation somehow captures the obstruction to 
commutativity of this diagram.

One special case of this theorem is when $F$ is a real quadratic field, and $\pi$ corresponds 
to a Hilbert modular form of CM type, then our period relations are formally the same period 
relations proved by Murty and Ramakrishnan in \cite{kumar-dinakar}. The reader should also 
compare our Theorem~\ref{thm:twisted} with the conjectures of Blasius and Panchiskin on the 
behaviour of Deligne's periods attached to a motive upon twisting it by Artin motives. See 
\S\ref{sec:remarks}.

\smallskip

A second aim of this paper is to consider Deligne's conjectures about the special values 
of the symmetric power $L$-functions attached to a holomorphic modular form. In 
\S\ref{sec:deligne} we recall the precise statement of the conjecture. See 
Conjecture~\ref{conj:deligne}. (See also our previous paper \cite{raghuram-shahidi} concerning 
the implications of some recent progress in Langlands functoriality for the special values of 
symmetric power $L$-functions.) If the modular cusp form is of CM type, namely, if the 
associated cuspidal representation is induced from a character $\chi$ of an imaginary quadratic 
extension, then we say the modular form is of dihedral type, and denote it by $\varphi_{\chi}$. 
Now Deligne's conjecture on special values of the symmetric power $L$-functions attached to 
dihedral forms is known via motivic techniques; this is because Deligne's main conjecture 
\cite[Conjecture 2.8]{deligne} is known if one considers only the motives as those attached to 
abelian varieties and the category used is that defined by using absolute Hodge cycles for 
morphisms. In \S\ref{sec:dihedral} we give an elementary proof in the dihedral case using only 
$L$-function techniques. In Theorem~\ref{thm:period-relations}, we prove 
another period relation, which after some formal functorial calculations, implies Deligne's
conjecture. We now proceed to describe these period relations. 

Given a holomorphic modular cusp form $\varphi$ on the upper half plane, let $u^{\pm}(\varphi)$
be the periods attached to $\varphi$ by Shimura \cite{shimura2}. The critical values
of $L(s,\varphi)$ are described in terms of $u^{\pm}(\varphi)$ \cite[Theorem 1]{shimura2}. 
Now, if $\varphi = \varphi_{\chi}$ is dihedral, then one can check that the $r$-th symmetric 
power
$L$-function attached to $\varphi_{\chi}$ is essentially a product of $L$-functions of 
$\varphi_{\chi^n}$ for various powers of $\chi$.  
See Corollary~\ref{cor:dihedral} for the precise factorization. 
One can deduce Deligne's conjecture if one relates the periods of $\varphi_{\chi^n}$ to 
the periods of $\varphi_{\chi}$. This is the content of our second main theorem of this paper. 
In \S\ref{sec:dihedral} we prove the following

\smallskip
\noindent{\bf Theorem~\ref{thm:period-relations}.} {\it 
Let $\chi$ be a Hecke character of $K$, where $K$ is
an imaginary quadratic extension of ${\mathbb Q}$. Assume that $\chi_{\infty}(z) = 
(z/|z|)^{k-1}$ for an integer $k \geq 2$. Let $\varphi_{\chi}$ be the corresponding 
primitive modular cusp form. Then we have 
$$
u^+(\varphi_{\chi^n}) \sim_{{\mathbb Q}(\chi)} u^+(\varphi_{\chi})^n,\ 
{\rm and}\ \  
u^-(\varphi_{\chi^n}) \sim_{{\mathbb Q}(\chi)} u^+(\varphi_{\chi})^n \gamma(\omega_K)
$$
where $\gamma(\omega_K)$ is the Gauss sum
of the quadratic Hecke character $\omega_K$ of ${\mathbb Q}$ attached to $K$ by class field 
theory; and $\sim_{{\mathbb Q}(\chi)}$ means up to an element of ${\mathbb Q}(\chi)$.} 

\smallskip

Again, we prove a stronger ${\rm Aut}({\mathbb C})$-equivariant version of this 
period relation. The proof of this theorem is by induction on $n$, while using 
the Rankin-Selberg $L$-function attached to $\varphi_{\chi^n} \times \varphi_{\chi}$; the proof 
also uses some
well known nonvanishing results for twists of $L$-functions (see Lemma~\ref{lem:nonvanishing}).

\smallskip

In \S\ref{sec:sym4} we take up the theme of our paper \cite{raghuram-shahidi}, of using 
Langlands functoriality to special values of $L$-functions, as applied to the case of the 
symmetric fourth lifting of a holomorphic modular form. We know from the work of Kim \cite{kim} 
that given a cuspidal representation $\pi$ of ${\rm GL}_2$, ${\rm Sym}^4(\pi)$ exists as an 
automorphic representation of ${\rm GL}_5$. The hope is to be able to use this, in conjunction 
with the recent work of Mahnkopf \cite{mahnkopf2} on the special values of standard 
$L$-functions on ${\rm GL}_n$ over ${\mathbb Q}$, to prove Deligne's conjecture on the special 
values of symmetric fourth power $L$-function attached to a modular form. Pursuing this 
line of 
thought, in \S\ref{sec:sym4-specialvalues}, we have recorded the current status of what is 
known, and what are some of the impediments; which in turn may be construed as an impetus for 
future work.

It might help the reader to know that \S\ref{sec:periods} and \S\ref{sec:dc} are quite 
independent of each other. However, the relatively short \S\ref{sec:sym4} depends on 
both \S\ref{sec:periods} and \S\ref{sec:dc}. 

\medskip

{\SMALL
\noindent {\it Acknowledgements:} We are grateful to Laurent Clozel, Paul Garrett, Joachim Mahnkopf, Dipendra~Prasad,  
Dinakar~Ramakrishnan and David Vogan for helpful correspondence. The first author thanks 
the warm hospitality of Purdue University. Both authors would like to thank Steve Kudla, 
Michael Rapoport 
and Joachim Schwermer for the invitation to spend some time in the stimulating atmosphere of 
the Erwin Schr\"odinger Institute in Vienna, where the work took its final form. This work is 
partially supported by the Vaughn foundation for (A.R.), and by NSF grants DMS-0200325 
and DMS-0700280 for (F.S.).}

\section{Periods of cusp forms}
\label{sec:periods}

\subsection{Notation and some preliminaries}
\label{sec:prelims}

For a number field $F$, we let ${\mathbb A}_F$ stand for its ad\`ele ring, and
${\mathbb I}_F = {\mathbb A}_F^{\times}$ be its group of id\`eles. 
We let $|\!|\ |\!|_F : {\mathbb I}_F \to {\mathbb R}_{>0}$ be the ad\`elic norm
defined by $|\!|x|\!|_F = \prod_v |x_v|_v$, with $v$ running over all places of $F$,
and the local absolute values all being the normalized ones. 
When there is no confusion about the base field $F$, we will drop the 
subscript $F$ from ${\mathbb A}_F$, ${\mathbb I}_F$, and $|\!|\ |\!|_F$. For any finite 
set $S$ of places
of $F$ we
use a superscript $S$ to denote a product outside $S$, and a subscript $S$ to
denote a product within $S$. For example, if $S_{\infty}$ stands for the
set of all infinite places of $F$, then the ring of finite ad\`eles
is ${\mathbb A}_F^{S_{\infty}}$ and will be denoted
${\mathbb A}_{F,f}$ or simply as ${\mathbb A}_f$.
We let $S_r$ stand for the set of real places and so $S_c := S_{\infty}-S_r$ is
the set of complex places. Let $r_1$ (respectively
$r_2$) denote the number of real (respectively complex) places of $F$;
the degree of $F$ is $d_F := [F:{\mathbb Q}] = r_1+2r_2.$ 

Let $G = {\rm GL}_n$, and let $Z = Z_n$ be the center of $G$, both regarded as $F$-groups. 
Let $G_{\infty} = G(F \otimes {\mathbb R}) = G({\mathbb R})^{r_1} \times
G({\mathbb C})^{r_2}$. 
Following Borel--Jacquet \cite[\S4.6]{borel-jacquet}, we say an irreducible 
representation of $G({\mathbb A})$ is automorphic if it is isomorphic to an 
irreducible subquotient of the representation of $G({\mathbb A})$ on its
space of automorphic forms. We say an automorphic representation is cuspidal 
if it is a subrepresentation of the representation of $G({\mathbb A})$ on 
the space of cusp forms $\mathcal{A}_{\rm cusp}(G(F)\backslash G({\mathbb A}))$.
For an automorphic representation $\pi$ of $G({\mathbb A})$, we have
$\pi = \pi_{\infty} \otimes \pi_f$, where $\pi_{\infty} =
\otimes_{v \in S_{\infty}} \pi_v$ is a representation of $G_{\infty}$,
and $\pi_f = \otimes_{v \notin S_{\infty}} \pi_v$
is a representation of $G({\mathbb A}_f)$.

By a Hecke character $\xi$ of $F$, we mean a continuous unitary character of
$F^*\backslash {\mathbb I}_F$. We follow the terminology as in Neukirch's book
\cite[\S VII.6]{neukirch}. Such a character admits a module of definition, say 
$\mathfrak{m}$, which is an integral ideal of $F$. 
If $\xi$ is a Hecke character modulo $\mathfrak{m}$, then we will
also identify $\xi$ with the corresponding Gr\"o\ss{}encharakter modulo $\mathfrak{m}$ as in
\cite[VII.6.14]{neukirch}.

\subsection{Definition of the periods}
\label{sec:periods-defn}

Let $F$ be a number field. The purpose of this section is to define certain periods 
attached to a regular algebraic cuspidal automorphic representation $\pi$ of ${\rm 
GL}_n({\mathbb A}_F)$. This definition is due to Harder \cite{harder} for ${\rm GL}_2$, 
and is due to Mahnkopf \cite{mahnkopf2} in the case $F = {\mathbb Q}$. 
(We refer the reader to Clozel \cite{clozel2} for the definitions of a cuspidal
representation being regular and algebraic.) 

Before we get into the details of the definition, we very roughly indicate the 
ingredients needed in making the definition. We will have a number field $E$. We will 
have two ${\mathbb C}$-vector spaces $V_1$ and $V_2$ with $E$-structures $V_1^0$ and 
$V_2^0$ respectively. (By $V_i^0$ being an $E$-structure for $V_i$, we mean an 
$E$-subspace such that the canonical map $V_i^0 \otimes_E {\mathbb C} \to V_i$ is an 
isomorphism.) In our situation, the spaces $V_i$ will be representation spaces, and 
not merely vector spaces, and the $E$-structures will be unique up to homotheties. 
Finally, we will have a comparison isomorphism $\phi : V_1 \to V_2$. The period attached 
to $\phi$, denoted $p(\phi)$, is a nonzero complex number such that $\phi(V_1^0) = 
p(\phi)V_2^0$. Observe that $p(\phi)$ is a well defined element in ${\mathbb C}^*/E^*$. 
For us, the number field $E$ will be the rationality field of $\pi$, the space 
$V_1$ will be the Whittaker model of $\pi$, and the space $V_2$ will be a certain 
cohomology space (to which $\pi$ will have nonzero contribution), and the comparison
isomorphism $\phi$ will be related to taking the Fourier coefficient of a cusp form
in the space of $\pi$. We now proceed to make all this precise.

The first ingredient we need is the {\it rationality field} of 
$\pi$, or really, $\pi_f$. The definitive reference is Clozel \cite[Chapter 
3]{clozel2}. Given $\pi$, suppose $V$ is the representation space of $\pi_f$, any 
$\sigma \in {\rm Aut}({\mathbb C})$ defines a representation $\pi^{\sigma}_f$ 
on $V \otimes_{\mathbb C} {\mathbb C}_{\sigma^{-1}}$ where $G({\mathbb A}_f)$ 
acts on the first factor. Let $\mathcal{S}(\pi_f)$ be the subgroup of ${\rm 
Aut}({\mathbb C})$ consisting of all $\sigma$ such that $\pi^{\sigma}_f \simeq 
\pi_f$. Define the rationality field ${\mathbb Q}(\pi_f)$ of $\pi_f$ as the 
subfield of ${\mathbb C}$ fixed by $\mathcal{S}(\pi_f)$; we denote this as 
${\mathbb Q}(\pi_f) = {\mathbb C}^{\mathcal{S}(\pi_f)}$. For example, if $\chi$ 
is a Dirichlet character, also thought of as an id\`ele class character, then 
${\mathbb Q}(\chi_f)$ is the field ${\mathbb Q}(\{\mbox{Values of $\chi$}\})$. 
Similarly, if $\varphi$ is a primitive holomorphic cusp form on the upper half 
plane, of even weight $2k$, for the Hecke congruence subgroup $\Gamma_0(N)$, with 
Fourier expansion $\varphi(z) = \sum_{n=1}^{\infty} a_nq^n$, and if $\pi = 
\pi(\varphi)$ is the cuspidal automorphic representation associated to $\varphi$, 
then ${\mathbb Q}(\pi_f) = {\mathbb Q}(\{a_n : n \geq 1\})$--the field generated 
by all the Fourier coefficients of $\varphi$. (See \cite{waldspurger1}.) In this 
example, the weight is assumed to be even to ensure that $\pi$ is algebraic. If 
the weight is odd, the same is true with $\pi$ replaced by 
$\pi \otimes |\!| \ |\!|^{-1/2}$. The main results that we need about the rationality 
field is stated in the following theorem. (See \cite[Th\'eor\`eme 3.13]{clozel2} and 
\cite[Chapter I]{waldspurger1}.)

\begin{thm}[Eichler, Shimura, Harder, Waldspurger, Clozel]
Let $\pi$ be a regular algebraic cuspidal automorphic representation of 
${\rm GL}_n({\mathbb A}_F)$. Then, 
\begin{enumerate}
\item ${\mathbb Q}(\pi_f)$ is a number field. 
\item $\pi_f$ admits a ${\mathbb Q}(\pi_f)$-structure, which is unique up to 
homotheties. 
\item For any $\sigma \in {\rm Aut}({\mathbb C})$, $\pi^{\sigma}_f$ is the finite part
of a cuspidal automorphic representation (which we denote by $\pi^{\sigma}$). 
\end{enumerate}
\end{thm}

The next ingredient we need is {\it the Whittaker model of $\pi_f$ and a 
semilinear action of ${\rm Aut}({\mathbb C})$ on this space}, which will commute 
with the action of ${\rm GL}_n({\mathbb A}_f)$.
Toward this, we fix a nontrivial character $\psi$ of $F\backslash {\mathbb A}_F$. (For 
convenience we fix $\psi$ as in Tate's thesis, namely,  
$\psi(x) = e^{2\pi i \Lambda(x)}$ with the $\Lambda$ as defined in 
\cite[\S 4.1]{tate}.) We can write $\psi = \psi_{\infty}\otimes \psi_f$ 
(the meaning 
and notation being the obvious one). We let $W(\pi, \psi)$ be the Whittaker model of 
$\pi$, and this factors as $W(\pi,\psi) = W(\pi_{\infty},\psi_{\infty}) \otimes 
W(\pi_f,\psi_f)$.
There is a semilinear action of ${\rm Aut}({\mathbb C})$ on $W(\pi_f,\psi_f)$ 
which is defined as follows. 
(See \cite[pp. 79-80]{harder} or \cite[pp. 594]{mahnkopf2}.) 
That the values of $\psi$ are all roots  
of unity suggests that we consider the cyclotomic character
{\small
$$
\begin{array}{llllllclc}
{\rm Aut}({\mathbb C}/{\mathbb Q}) & \to & 
{\rm Gal}(\overline{\mathbb Q}/{\mathbb Q}) & \to&
{\rm Gal}({\mathbb Q}(\mu_{\infty})/{\mathbb Q}) & \to &
\widehat{{\mathbb Z}}^{\times} \simeq \prod_p {\mathbb Z}_p^{\times} & \subset & 
\prod_p \prod_{\mathfrak{p}|p} \mathcal{O}_{\mathfrak{p}}^{\times} \\
\sigma & \mapsto & \sigma |_{\overline{\mathbb Q}} & \mapsto & 
\sigma |_{{\mathbb Q}(\mu_{\infty})} & \mapsto & t_{\sigma} & \mapsto & t_{\sigma}
\end{array}
$$}
where the last inclusion is the one induced by the diagonal embedding of 
${\mathbb Z}_p^{\times}$ into $\prod_{\mathfrak{p}|p} 
\mathcal{O}_{\mathfrak{p}}^{\times}$. (Here $\mathfrak{p}$ is a prime of $F$ above $p$,
and $\mathcal{O}_{\mathfrak{p}}$ is the ring of integers of the completion 
$F_{\mathfrak{p}}$ of $F$ at $\mathfrak{p}$.) The element $t_{\sigma}$ at the end can
be thought of as an element of ${\mathbb A}_f^{\times} = {\mathbb I}_f$. 
Let $t_{\sigma,n}$ denote the diagonal matrix 
${\rm diag}(t_{\sigma}^{-(n-1)}, t_{\sigma}^{-(n-2)},\dots,1)$ regarded as an 
element of  ${\rm GL}_n({\mathbb A}_f)$. 
For $\sigma \in {\rm Aut}({\mathbb C})$ and 
$\phi \in W(\pi_f,\psi_f)$, define the function $W_{\sigma}(\phi)$ by
$$
W_{\sigma}(\phi)(g_f) = \sigma(\phi(t_{\sigma,n}g_f))
$$
for all $g_f \in {\rm GL}_n({\mathbb A}_f)$. It is easily seen that
$W_{\sigma}$ is a $\sigma$-linear ${\rm GL}_n({\mathbb A}_f)$-equivariant
isomorphism from $W(\pi_f,\psi_f)$ onto $W(\pi^{\sigma}_f, \psi_f)$. 
For any finite extension $E/{\mathbb Q}(\pi_f)$ we have
an $E$-structure on $W(\pi_f,\psi_f)$ by taking invariants: 
$$
W(\pi_f,\psi_f)_E = W(\pi_f,\psi_f)^{{\rm Aut}({\mathbb C}/E)}.
$$

As a matter of notation, given a ${\mathbb C}$-vector space $V$,
and given a subfield $E \subset {\mathbb C}$, we will let $V_E$
stand for an $E$-structure on $V$ (if there is one at hand). 
Fixing an $E$-structure
gives an action of ${\rm Aut}({\mathbb C}/E)$ on $V$, 
by making it act on the 
second factor in $V = V_E \otimes_E {\mathbb C}$. Having fixed an 
$E$-structure, for any extension $E'/E$, we have a canonical $E'$-structure
by letting $V_{E'} = V_E \otimes_E E'$.  
Further, as a notational convenience, when we talk of Whittaker
models, we will henceforth suppress the additive character $\psi$, since
that has been fixed once and for all; for example, $W(\pi_f)$ will denote 
$W(\pi_f,\psi_f)$. Also, we will denote the map $W_{\sigma}$ 
simply by $\sigma$.

As mentioned earlier, the periods come via a 
comparison of $W(\pi_f)_E$ with a rational structure on a 
suitable cohomology space. We now describe this cohomology space. 
Recall that $G = G_n = {\rm GL}_n$ and the center of $G$ is denoted $Z_n$ or $Z$. 
Let $\mathfrak{g}_{\infty}$ be the Lie algebra of $G_{\infty}.$ 
Let $K_{\infty} = \otimes_{v \in S_{\infty}} K_v$ 
where $K_v = Z_n({\mathbb R})O_n({\mathbb R})$ if $v$ is real, and 
$K_v = Z_n({\mathbb C})U_n({\mathbb C})$ if $v$ is complex. 
Let $K^0 = K_{\infty}^0$ be the topological connected component of $K_{\infty}.$ 
Note that $K_{\infty}/K_{\infty}^0 \simeq ({\mathbb Z}/2{\mathbb Z})^{r_1}$.
Let $b_n^{\mathbb R}$ be $n^2/4$ if $n$ is even, and $(n^2-1)/4$ if $n$ is odd.
We also let $b_n^{\mathbb C}$ be $n(n-1)/2$. Now we define
$b = r_1b_n^{\mathbb R} + r_2b_n^{\mathbb C}.$
The integer $b$ depends only on the base field $F$ and the rank $n$ of 
${\rm GL}_n$. It is the {\it bottom degree} of the so called {\it cuspidal range} 
for ${\rm GL}_n$ as an $F$-group. The next ingredient we need 
in defining the period is {\it relative Lie algebra 
cohomology of $\pi$ in degree $b$.} 
For a $(\mathfrak{g}_{\infty}, K_{\infty}^0)$-module $U$, we let 
$H^*(\mathfrak{g}_{\infty}, K_{\infty}^0; U)$ stand for relative Lie 
algebra cohomology
of $U$, for the definition and properties of which 
we refer the reader to Borel and Wallach's book \cite{borel-wallach}. 
Given a representation $\tau$ of $G_{\infty}$, by 
$H^*(\mathfrak{g}_{\infty}, K_{\infty}^0; \tau)$, we will mean the cohomology
of the $(\mathfrak{g}_{\infty}, K_{\infty}^0)$-module consisting of 
smooth $K_{\infty}$-finite vectors of $\tau$.

Let $T = T_n$ denote the maximal torus of ${\rm GL}_n$ consisting of diagonal 
matrices. We regard $T$ as an $F$-group, and let 
$T_{\infty} = T(F \otimes {\mathbb R}) = T({\mathbb R})^{r_1} \times 
T({\mathbb C})^{r_2}$. 
We let $B = B_n$ stand for the Borel subgroup of $G$ 
of upper triangular matrices. This defines $B_{\infty}$. We let $X(T_{\infty})$ 
stand for the group of all algebraic characters of $T_{\infty}$. We let 
$X^+(T_{\infty})$ stand for the subset of $X(T_{\infty})$ consisting of all those characters
which are dominant with respect to $B_{\infty}$.
A weight $\mu \in 
X^+(T_{\infty})$ may be described as follows: Let $\mu = (\mu_v)_{v \in S_{\infty}}$, with 
$\mu_v \in X(T_v)$. 
If $v \in S_r$, then $\mu_v = (p_1,\dots,p_n)$, $p_i \in {\mathbb Z}$, 
$p_1 \geq p_2 \geq \cdots \geq p_n$, 
and the character is: if $t = {\rm diag}(t_1,\dots,t_n) \in T(F_v) = T({\mathbb R})$,
then $t \mapsto \prod_i t_i^{p_i}.$ If $v \in S_c$, then let $\{\iota_v, \bar{\iota}_v\}$
be the corresponding complex embeddings of $F$. Identify $F_v$ with ${\mathbb C}$
via $\iota_v$ (say). In this case, $\mu_v$ is a pair of sequences $(\mu_{\iota_v}, 
\mu_{\bar{\iota}_v})$, with $\mu_{\iota_v} = (p_1,\dots,p_n)$, 
 $p_i \in {\mathbb Z}$, $p_1 \geq p_2 \geq \cdots \geq p_n$; likewise
$\mu_{\bar{\iota}_v} = (q_1,\dots,q_n)$ with similar conditions on the $q_i$'s; 
the character $\mu_v$ is: if $t = {\rm diag}(z_1,\dots,z_n) \in T(F_v) = T({\mathbb C}),$
then $t \mapsto \prod_i z_i^{p_i} \bar{z}_i^{q_i}$. (Here $\bar{z_i}$ is the complex 
conjugate of $z_i$.) 
For such a character $\mu$, we define a finite dimensional
representation $(\rho_{\mu},M_{\mu})$ of $G_{\infty}$ as follows.
For $v \in S_r$, let $(\rho_{\mu_v}, M_{\mu_v})$ 
be the irreducible representation of $G(F_v) = G({\mathbb R})$ with 
highest weight $\mu_v$. For $v \in S_c$, let $(\rho_{\mu_v}, M_{\mu_v})$ 
be the representation of $G(F_v) = G({\mathbb C})$ defined as 
$\rho_{\mu_v} = \rho_{\mu_{\iota_v}} \otimes \rho_{\mu_{\bar{\iota}_v}}$, where
$\rho_{\mu_{\iota_v}}$ is the irreducible representation with highest weight 
$\mu_{\iota_v}$, and similarly $\rho_{\mu_{\bar{\iota}_v}}$.
Now we let $\rho_{\mu} = \otimes_{v \in S_{\infty}} \rho_{\mu_v}.$
Since $\pi$ is a regular algebraic cuspidal automorphic
representation of ${\rm GL}_n({\mathbb A})$, we have 
from the proof of \cite[Th\'eor\`eme 3.13]{clozel2} that there is a dominant
algebraic weight $\mu \in X^+(T_{\infty})$ 
such that 
$
H^*(\mathfrak{g}_{\infty}, K_{\infty}^0;\pi_{\infty} \otimes \rho_{\mu}^{\vee})
\neq 0.$
In defining the periods, we will be looking at 
$H^b(\mathfrak{g}_{\infty}, K_{\infty}^0;\pi_{\infty} \otimes \rho_{\mu}^{\vee})$.

The group $K_{\infty}/K_{\infty}^0 \simeq ({\mathbb Z}/2{\mathbb Z})^{r_1}$ acts 
on $H^b(\mathfrak{g}_{\infty}, K_{\infty}^0;\pi_{\infty} \otimes M_{\mu}^{\vee})$.
We consider certain isotypic components for this action. Consider an 
$r_1$ tuple of signs indexed by the set $S_r$ of real places in $S_{\infty}$. 
Let $\epsilon = (\epsilon_v)_{v\in S_r} \in 
\{1\!\!1, {\rm sgn}\}^{r_1} = (K_{\infty}/K_{\infty}^0)\widehat{}$. 
If $n$ is even then there are no restrictions on $\epsilon$, however, if
$n$ is odd then $\pi$ uniquely determines an $\epsilon$, in that we 
let $\epsilon_v = \omega_{\pi_v}|_{\pm 1}\ \cdot {\rm sgn}^{{\rm wt}(\mu_v)/2}$, where
the weight ${\rm wt}(\mu_v)$ of $\mu_v$ is defined in \cite[(3.1)]{mahnkopf2}.
(If $n$ is odd, then for $v \in S_r$, $\epsilon_v$ is simply the central character
of $\pi_v\otimes\rho_{\mu_v}^{\vee}$ restricted to $\{\pm 1\}$, since the parity
of $n$ means that $K_v/K_v^0 \simeq \{\pm 1\} \subset Z_n(F_v).$)
Let $H^b(\mathfrak{g}_{\infty},K_{\infty}^0; \pi_{\infty}\otimes 
M_{\mu}^{\vee})(\epsilon)$ be the corresponding isotypic 
component. This isotypic component is one dimensional. This can be seen, by 
using \cite[(3.2)]{mahnkopf2} for the real places, \cite[Lemme 3.14]{clozel2} for
the (real and) complex places, and the K\"unneth 
theorem for Lie algebra cohomology \cite[I.1.3]{borel-wallach}, as follows:
\begin{eqnarray*}
& & H^b(\mathfrak{g}_{\infty},K_{\infty}^0; \pi_{\infty}\otimes
M_{\mu}^{\vee})(\epsilon) \\
& = &
\bigoplus_{\sum a_v = b}\left( 
\bigotimes_{v \in S_r}
H^{a_v}(\mathfrak{g}_v, K_v^0; \pi_v \otimes M_{\mu_v}^{\vee})(\epsilon_v)
\bigotimes_{v \in S_c}
H^{a_v}(\mathfrak{g}_v,K_v^0; \pi_v \otimes M_{\mu_v}^{\vee})
\right) \\
& = & 
\bigotimes_{v \in S_r}
 H^{b_n^{\mathbb R}}(\mathfrak{g}_v,K_v^0; \pi_v \otimes M_{\mu_v}^{\vee})
(\epsilon_v)
\bigotimes_{v \in S_c}
H^{b_n^{\mathbb C}}(\mathfrak{g}_v,K_v^0; \pi_v \otimes M_{\mu_v}^{\vee}).
\end{eqnarray*}
In the summation, 
only one term survives, because for all other summands, at least one of the 
$a_v$ has to be less than $b_n^{\mathbb R}$ or $b_n^{\mathbb C}$, and by 
\cite[Lemme 3.14]{clozel2} the corresponding factor vanishes. 
We fix a generator ${\bf w}_{\infty} = {\bf w}(\pi_{\infty}, \epsilon)$ for 
this one dimensional space $H^b(\mathfrak{g}_{\infty},K_{\infty}^0; 
\pi_{\infty}\otimes
M_{\mu}^{\vee})(\epsilon)$.

We have the following comparison isomorphism of the Whittaker model $W(\pi_f)$ with a
global version of the above cohomology space. We
let $\mathcal{F}_{\pi_f, \epsilon, {\bf w}_{\infty}}$ denote the compositum of
the three isomorphisms: 
\begin{eqnarray*}
W(\pi_f) & \longrightarrow & 
W(\pi_f) \otimes 
H^b(\mathfrak{g}_{\infty},K_{\infty}^0; W(\pi_{\infty}) \otimes 
M_{\mu}^{\vee})(\epsilon) \\
& \longrightarrow & 
H^b(\mathfrak{g}_{\infty},K_{\infty}^0; W(\pi) \otimes 
M_{\mu}^{\vee})(\epsilon) \\
& \longrightarrow & 
H^b(\mathfrak{g}_{\infty},K_{\infty}^0; V_{\pi} \otimes M_{\mu}^{\vee})(\epsilon),
\end{eqnarray*}
where the first map is $w_f \mapsto w_f \otimes {\bf w}_{\infty}$; 
the second map is the obvious one; and the third map is the map induced in cohomology 
by the inverse
of the map which gives the Fourier coefficient of a cusp form in $V_{\pi}$--the space
of functions in $\mathcal{A}_{{\rm cusp}}(G(F) \backslash G({\mathbb A}))$ which 
realizes $\pi$. 

We now describe a rational structure on 
$H^b(\mathfrak{g}_{\infty},K_{\infty}^0; V_{\pi} \otimes 
M_{\mu}^{\vee})(\epsilon)$, by relating relative Lie algebra cohomology with the cohomology
of locally symmetric spaces. (See \cite[pp.128--129]{clozel2}, \cite[\S 3.2]{mahnkopf2}.)
Let $K_f$ be any open compact subgroup of $G({\mathbb A}_f)$. Consider the manifold
$$
S(K_f) = G(F)\backslash G({\mathbb A})/K_{\infty}^0K_f.
$$
This is typically a finite disjoint union of manifolds like 
$\Gamma \backslash G_{\infty}/K_{\infty}^0$. 
For a dominant algebraic weight $\mu \in X^+(T_{\infty})$ the 
corresponding finite dimensional representation $(\rho_{\mu},M_{\mu})$ of $G_{\infty}$ 
is defined over ${\mathbb Q}$. Fix a ${\mathbb Q}$-structure $M_{\mu, {\mathbb Q}}$
on $M_{\mu}$ which gives a canonical $E$-structure 
for any extension $E/{\mathbb Q}$ by 
$M_{\mu, E} = M_{\mu, {\mathbb Q}} \otimes E$. 
Let $\mathcal{M}_{\mu, E}$ be the associated locally constant sheaf on $S(K_f)$. 
For brevity, we also let $M_{\mu} = M_{\mu, {\mathbb C}}$ and similarly 
$\mathcal{M}_{\mu} = \mathcal{M}_{\mu, {\mathbb C}}$. 
We consider the direct limit of various cohomology groups
$$
H^{\bullet}_?(\tilde{S},\mathcal{M}^{\vee}_{\mu}) = 
\underrightarrow{\rm lim}\ 
H^{\bullet}_?(S(K_f), \mathcal{M}^{\vee}_{\mu}),
$$
where the direct limit is taken over all open compact subgroups $K_f$ of $G({\mathbb 
A}_f)$,
and $? \in \{B, dR, c, {\rm cusp}\}$ meaning 
singular (Betti) cohomology, or de Rham cohomology, or cohomology 
with compact supports, or cuspidal cohomology, respectively. 
Cuspidal cohomology injects into cohomology with compact supports 
$H^{\bullet}_{\rm cusp}(\tilde{S},\mathcal{M}^{\vee}_{\mu}) 
\hookrightarrow 
H^{\bullet}_c(\tilde{S},\mathcal{M}^{\vee}_{\mu})$
and the latter is canonically a module for 
${\rm Aut}({\mathbb C}) \times G({\mathbb A}_f)$
with commuting actions for the different groups. The image of cuspidal
cohomology is defined over ${\mathbb Q}$. Also, cuspidal cohomology
decomposes into a direct sum over cuspidal cohomological representations
and a rational structure on any summand is obtained by intersecting that
summand with a rational structure of the ambient space. (As a general reference
for all of this, see Clozel \cite[\S 3.5]{clozel2}.)

More precisely, by definition of cuspidal cohomology, we have 
\begin{equation}
\label{eqn:cuspidal-cohomology}
H^b_{\rm cusp}(\tilde{S},\mathcal{M}^{\vee}_{\mu}) \simeq 
H^b(\mathfrak{g}_{\infty},K_{\infty}^0; \mathcal{A}_{\rm cusp}(G(F)\backslash G({\mathbb 
A}))
\otimes M_{\mu}^{\vee}).
\end{equation}
From the decomposition of $\mathcal{A}_{\rm cusp}(G(F)\backslash G({\mathbb A}))$ into 
cuspidal automorphic representations, we deduce that the right hand side decomposes into
a direct sum
$$
H^b(\mathfrak{g}_{\infty},K_{\infty}^0; \mathcal{A}_{\rm cusp}(G(F)\backslash G({\mathbb 
A}))
\otimes M_{\mu}^{\vee}) \simeq 
\bigoplus_{\pi \in {\rm Coh}(G, \mu^{\vee})}
H^b(\mathfrak{g}_{\infty},K_{\infty}^0; V_{\pi} \otimes M_{\mu}^{\vee}).
$$ 
This also defines the notation ${\rm Coh}(G, \mu^{\vee})$ as the set consisting of all $\pi$
having a nonzero contribution in the right hand side. 
Now consider the action of 
$K_{\infty}/K_{\infty}^0$, and further decompose each summand into its isotypic components.  
Let $\epsilon \in (K_{\infty}/K_{\infty}^0){}^{\widehat{}}$ be as before, i.e.,
if $n$ is even then $\epsilon$ can be any character, and 
if $n$ is odd, then $\epsilon$ is uniquely determined by $\pi$. 
We let $\tilde{\pi} = \pi_f \otimes \epsilon$, and 
denote the inverse image of $H^b(\mathfrak{g}_{\infty},K_{\infty}^0; V_{\pi} \otimes 
M_{\mu}^{\vee})(\epsilon)$ across the isomorphism in (\ref{eqn:cuspidal-cohomology}) by
$H^b_{\rm cusp}(\tilde{S},\mathcal{M}^{\vee}_{\mu})(\tilde{\pi})$. We have 
$$
H^b_{\rm cusp}(\tilde{S},\mathcal{M}^{\vee}_{\mu}) \simeq
\bigoplus_{\pi \in {\rm Coh}(G, \mu^{\vee})} 
\bigoplus_{\epsilon}
H^b_{\rm cusp}(\tilde{S},\mathcal{M}^{\vee}_{\mu})(\tilde{\pi}),
$$
where in the second summation $\epsilon$ is as described above. (See also
\cite[(3.13)]{mahnkopf2}.)

We now have the following description of the rational
structures. The Betti cohomology spaces $H^b_B(\tilde{S}, \mathcal{M}^{\vee}_{\mu})$ 
are naturally defined over ${\mathbb Q}$, since the coefficient system 
admits a ${\mathbb Q}$-structure \cite[p.122]{clozel2}. 
(This will be exploited in the proof of Proposition~\ref{prop:automorphic-rational}.) 
The ${\mathbb Q}$-structure on Betti cohomology induces a ${\mathbb Q}$-structure on 
cohomology with compact support \cite[exact triangle on p.123]{clozel2} which we denote
by $H^b_c(\tilde{S}, \mathcal{M}^{\vee}_{\mu})_{\mathbb Q}$.
This in turn induces a 
${\mathbb Q}$-structure on the full space of cuspidal cohomology
(\cite[Th\'eor\`em 3.19]{clozel2})
\begin{equation}
\label{eqn:rational-fullcusp}
H^b_{\rm cusp}(\tilde{S},\mathcal{M}^{\vee}_{\mu})_{\mathbb Q} = 
H^b_{\rm cusp}(\tilde{S},\mathcal{M}^{\vee}_{\mu}) \cap 
H^b_c(\tilde{S}, \mathcal{M}^{\vee}_{\mu})_{\mathbb Q}. 
\end{equation}
We get for each summand of cuspidal cohomology \cite[Lemme 3.2.1]{clozel2}
\begin{equation}
\label{eqn:rational-summand}
H^b_{\rm cusp}(\tilde{S},\mathcal{M}^{\vee}_{\mu})(\tilde{\pi})_E = 
H^b_{\rm cusp}(\tilde{S},\mathcal{M}^{\vee}_{\mu})(\tilde{\pi}) \cap 
H^b_{\rm cusp}(\tilde{S}, \mathcal{M}^{\vee}_{\mu})_E 
\end{equation}
for any extension $E/{\mathbb Q}(\pi_f)$. 
We can transport the rational structures (\ref{eqn:rational-fullcusp}) and 
(\ref{eqn:rational-summand}) across the 
identifications with relative Lie algebra cohomology 
to get rational structures on the latter spaces as: 
$$
H^b(\mathfrak{g}_{\infty}, K_{\infty}^0; \mathcal{A}_{\rm cusp}
(G(F)\backslash G({\mathbb A}))
\otimes M_{\mu}^{\vee})_{\mathbb Q} :\simeq
H^b_{\rm cusp}(\tilde{S},\mathcal{M}^{\vee}_{\mu})_{\mathbb Q}, 
$$
and for any extension $E/{\mathbb Q}(\pi_f)$ we have
\begin{eqnarray*}
H^b(\mathfrak{g}_{\infty}, K_{\infty}^0; V_{\pi} \otimes M_{\mu}^{\vee}) 
(\epsilon)_E & = & 
H^b(\mathfrak{g}_{\infty},K_{\infty}^0; V_{\pi} \otimes M_{\mu}^{\vee}) 
(\epsilon) \cap \\ & &
H^b(\mathfrak{g}_{\infty}, K_{\infty}^0; \mathcal{A}_{\rm cusp}(G(F)\backslash 
G({\mathbb A})) \otimes M_{\mu}^{\vee})_E.
\end{eqnarray*}

We recall that ${\rm Aut}({\mathbb C})$ acts on objects indexed by $S_{\infty}$ by 
acting on the embeddings of $F$ into ${\mathbb C}$. For the precise definition of 
the action of ${\rm Aut}({\mathbb C})$ on $\pi_{\infty}$ see 
\cite[\S 3.3]{clozel2}; this then defines an action of ${\rm Aut}({\mathbb C})$ on
$\mu$ and $\epsilon$ (they are determined by $\pi$), and finally, for 
$\sigma \in {\rm Aut}({\mathbb C})$ we define ${\bf w}_{\infty}^{\sigma}$ as 
${\bf w}(\pi_{\infty},\epsilon)^{\sigma} := 
{\bf w}(\pi_{\infty}^{\sigma},\epsilon^{\sigma})$. We are now in a position
to define the periods attached to $\pi$.

\begin{defn}
\label{defn:period}
Let $\pi$ be a regular algebraic cuspidal automorphic representation of 
${\rm GL}_n({\mathbb A}_F)$. Let $\mu \in X^+(T_{\infty})$ be such 
that $\pi \in {\rm Coh}(G, \mu^{\vee})$. Let $\epsilon$ be a character of 
$K_{\infty}/K_{\infty}^0$. (If $n$ is even then $\epsilon$ is any character, and if 
$n$ is odd, then $\pi$ uniquely determines $\epsilon$.) 
Let ${\bf w}_{\infty}$ be a generator of the one dimensional vector space 
$H^b(\mathfrak{g}_{\infty},K_{\infty}^0, \pi_{\infty} \otimes M_{\mu}^{\vee})(\epsilon)$. 
To this data $(\pi, \epsilon, {\bf w}_{\infty})$ is attached a period, that we denote
$p^{\epsilon}(\pi_f, {\bf w}_{\infty})$, which is a nonzero complex number such that the 
normalized 
map
$$
\mathcal{F}^0_{\pi_f, \epsilon, {\bf w}_{\infty}} := 
p^{\epsilon}(\pi_f, {\bf w}_{\infty})^{-1}
\mathcal{F}_{\pi_f, \epsilon, {\bf w}_{\infty}}
$$
is ${\rm Aut}({\mathbb C})$-equivariant: 
$$
\xymatrix{
W(\pi_f) \ar[rrr]^-{\mathcal{F}^0_{\pi_f, \epsilon, {\bf w}_{\infty}}}\ar[d]_{\sigma} & & &
H^b(\mathfrak{g}_{\infty},K_{\infty}^0; V_{\pi}\otimes M_{\mu}^{\vee})(\epsilon)
\ar[d]^{\sigma} \\
W(\pi^{\sigma}_f) 
\ar[rrr]^-{\mathcal{F}^0_{\pi^{\sigma}_f, \epsilon^{\sigma}, {\bf w}_{\infty}^{\sigma}}} 
& & &
H^b(\mathfrak{g}_{\infty},K_{\infty}^0; V_{\pi^{\sigma}}\otimes 
M_{\mu^{\sigma}}^{\vee})
(\epsilon^{\sigma})
}
$$
The period $p^{\epsilon}(\pi_f, {\bf w}_{\infty})$ is well defined only up to multiplication 
by elements of ${\mathbb Q}(\pi_f)^*$. 
\end{defn}

It is helpful to simplify our notation a bit. We begin by fixing generators for all the 
possible one dimensional relative Lie algebra cohomology spaces for the groups ${\rm 
GL}_n({\mathbb R})$ and ${\rm GL}_n({\mathbb C})$. Having done so, we have therefore fixed 
generators for the cohomology spaces for the group $G_{\infty}$. We also ask that this 
choice be compatible with twisting $\pi_{\infty}$ by algebraic unitary characters 
$\xi_{\infty}$ of $G_{\infty}$; this condition although crucial in the proof of 
Proposition~\ref{prop:diagram}, is not a serious constraint. (For example, for $G = {\rm 
GL}_2$, the choice of a generator for $H^1$ as in Waldspurger 
\cite[p.130-131]{waldspurger1} is already invariant under twisting.) This choice is 
compatible with the action of ${\rm Aut}({\mathbb C})$ on automorphic representations at 
infinity. Henceforth, we abbreviate $\mathcal{F}_{\pi_f, \epsilon, {\bf w}_{\infty}}$ by 
$\mathcal{F}_{\pi_f, \epsilon}$, similarly for the normalized maps, as well as 
$p^{\epsilon}(\pi_f, {\bf w}_{\infty})$ by $p^{\epsilon}(\pi_f)$, while keeping in mind 
that ${\bf w}_{\infty}$ has been chosen already. (For example, in the classical setting of 
modular forms, this is equivalent to fixing a generator ${\bf w}_{\infty}$ for $H^1$ of 
the discrete series representation of ${\rm GL}_2({\mathbb R})$ of lowest weight $k$, and 
now for all weight $k$ modular forms we work with this choice of ${\bf w}_{\infty}$.)

In terms of the un-normalized maps, we can describe the above commutative diagram by 
\begin{equation}
\label{eqn:unnormalized}
\sigma \circ \mathcal{F}_{\pi_f, \epsilon} = 
\left(\frac{\sigma(p^{\epsilon}(\pi_f))}
{p^{\epsilon^{\sigma}}(\pi_f^{\sigma})} \right) 
\mathcal{F}_{\pi_f^{\sigma}, \epsilon^{\sigma}} \circ \sigma.
\end{equation}

The periods $p^{\epsilon}(\pi_f)$ are expected to be intimately related to the special 
values of the $L$-function $L_f(s, \pi)$ attached to $\pi_f$; indeed, this is one of 
the main motivations for this paper. 
Making this relation explicit is part of a future work which we hope to report on in 
another paper. If $F = {\mathbb Q}$,
this is the main thrust of the work of Mahnkopf \cite{mahnkopf2}. Roughly speaking,
the quantity $\Omega_{\epsilon}(\pi, \pi', \eta)$ that shows up in 
\cite[Theorem 5.4]{mahnkopf2}--the main theorem of that paper--is given by 
$$
\Omega_{\epsilon}(\pi, \pi', \eta) = \frac{p^{\epsilon}(\pi_f)p^{\epsilon'}(\pi'_f)}
{L(0, \pi \times (\pi'\otimes\eta))}.
$$
For the precise definition of $\Omega_{\epsilon}(\pi, \pi', \eta)$,
which is considerably more involved and quite delicate, see \cite[(5.11)]{mahnkopf2}.
(In \S\ref{sec:sym4-specialvalues} we discuss a little more about certain issues posed by this 
auxiliary character $\eta$.)

Before we end this subsection we note that there is another approach to get an 
$E$-structure on Whittaker models. By using Waldspurger \cite[Lemme I.1]{waldspurger1} we 
get that there is an $E$-structure generated by all the local new vectors. It is an 
interesting exercise to compare the two different $E$-structures on $W(\pi_f)$ because 
that gives a {\it period} associated to $\pi_f$ (which is quite different from the Harder 
or Mahnkopf type periods). Whether this period has anything to do with special values of 
the $L$-function associated to $\pi$ is not clear at the moment.

\subsection{Behaviour under twisting}
\label{sec:period-twist}

The motivation for this section comes from our earlier work
\cite{raghuram-shahidi}. We formulated a conjecture which 
relates the special values of $L(s, {\rm Sym}^n \varphi, \xi)$--the symmetric
power $L$-function of a holomorphic cusp form $\varphi$ twisted by a Dirichlet character 
$\xi$--to the special values of $L(s, {\rm Sym}^n \varphi)$. We predict therein that a certain 
explicit power of the Gauss sum of $\xi$ relates the two special values. 
(See \cite[Conjecture 7.1]{raghuram-shahidi}.) 
Assuming Langlands' principle of  functoriality, this conjecture would {\it follow}
if we can prove a similar statement on the relation between the special values
of $L(s, \pi\otimes \xi)$, for a cuspidal representation $\pi$ of 
${\rm GL}_n({\mathbb A}_{\mathbb Q})$, and 
the special values of $L(s,\pi)$ itself. In view of the work of Mahnkopf, this translates to
proving a similar relation between the periods $P(\widetilde{\pi\otimes\xi})$ and 
$P(\tilde{\pi})$ associated by Mahnkopf to $\pi\otimes\xi$ and $\pi$ respectively. 
This relation would roughly take the form 
$$
P(\widetilde{\pi\otimes\xi}) \sim P(\tilde{\pi})
$$
where $\sim$ means up to some algebraic quantities over which one hopes to have as much
control as possible. 

The main aim of this section is to prove such a result for the periods $p^{\epsilon}(\pi_f)$ 
that we defined in the previous section. Albeit this statement is not trivial to prove, it 
should not be philosophically surprising, because, it seems unlikely to be able to introduce 
new transcendental quantities if one only allows twisting by `algebraic' characters. In a 
classical context, one should view this as that the transcendental part of the special value 
$L(2m,\chi)$, for an even Dirichlet character $\chi$ and an integer $m \geq 1$, is already 
captured by the special value $\zeta(2m)$ of the Riemann zeta function; and the quotient of 
the two special values is basically the Gauss sum of $\chi$. See \cite[\S 7]{raghuram-shahidi} 
for a discussion of such and other classical examples.

Before we state and prove the main theorem, we need some preliminaries on Hecke characters
and their Gauss sums. We keep all the notation introduced so far. Let $\xi$ be a 
Hecke character of $F$, and 
let $\xi^0$ be 
the associated primitive Hecke character and let the conductor of $\xi^0$ (and hence of $\xi$) 
be $\mathfrak{c}$. Let $\mathcal{D}_F$ be the absolute different of $F$, and let $y \in 
\mathfrak{c}^{-1}\mathcal{D}_F^{-1}$. Define the Gauss sum of $\xi$ by 
$$
\gamma(\xi_f,y) := \gamma(\xi^0_f,y) = 
\sum_{x \in (\mathcal{O}_F/\mathfrak{c})^{\times}} \xi_f(x)
e^{2 \pi i \, T_{F/{\mathbb Q}}(xy)}.
$$
In the right hand side, $\xi_f$ is the finite part of the Gr\"o\ss{}encharakter $\xi$.
This definition is slightly different from \cite[VII.6.3]{neukirch} and is more like the
definition of Gauss sum used by Shimura \cite[p.784]{shimura1}.
Observe that $\gamma(\xi_f, y)$ depends only on the class of $y$ in 
$\mathfrak{c}^{-1}\mathcal{D}_F^{-1}/\mathcal{D}_F^{-1}.$  Given a $\xi$, we will arbitrarily 
pick such an element $y$ for which $\gamma(\xi_f, y) \neq 0.$ Having chosen $y$ for $\xi$,
we will work with the same $y$ for every character of the form $\xi^{\sigma}$, where  
$\sigma \in {\rm Aut}({\mathbb C})$. This choice will not affect us in any serious way, 
because we will really be concerned with certain quotients involving Gauss sums, and 
such quotients will not depend on $y$. (See Lemma~\ref{lem:sigma-gauss} below.) In the 
notation
we will therefore suppress the dependence on $y$, and denote {\it the Gauss sum} of 
$\xi$ by $\gamma(\xi_f)$. 

Given a Hecke character $\xi$, we define its signature $\epsilon_{\xi} = (\epsilon_{\xi, 
v})_{v \in 
S_r} \in \{\pm 1\}^{r_1}$ as follows. 
For $v \in S_r$, define $\epsilon_{\xi, v} = \xi_v(-1)$. We will think of $\epsilon_{\xi}$
as a character of $K_{\infty}/K_{\infty}^0$. 
We can now state and prove one of the main 
results of this paper.

\begin{thm}
\label{thm:twisted}
Let $F$ be a number field and 
$\pi$ be a regular algebraic cuspidal automorphic representation of ${\rm GL}_n({\mathbb 
A}_F)$. 
Let $\mu$ be a dominant algebraic character of $T_{\infty}$ such 
that $\pi \in Coh(G,\mu^{\vee})$. Let $\epsilon$ be a character of 
$K_{\infty}/K_{\infty}^0$ as in \S\ref{sec:periods},
and let $p^{\epsilon}(\pi_f)$ be the period as in Definition~\ref{defn:period}. 
Let $\xi$ be an algebraic Hecke character of $F$ with signature $\epsilon_{\xi}$.
Then $\pi\otimes \xi$ is also a regular algebraic cuspidal automorphic representation, and 
the signature $\epsilon \cdot \epsilon_{\xi}:= \epsilon \otimes \epsilon_{\xi} $ is a character 
of $K_{\infty}/K_{\infty}^0$ that is permissible for 
$\pi\otimes \xi$ (which is an issue only when $n$ is odd), hence 
the period $p^{\epsilon\cdot \epsilon_{\xi}}(\pi_f\otimes\xi_f)$ is defined. We have 
the following relations:

\begin{enumerate}
\item For any $\sigma \in {\rm Aut}({\mathbb C})$ we have 
$$
\sigma\left(\frac
{p^{\epsilon \cdot \epsilon_{\xi}}(\pi_f\otimes\xi_f)}
{\gamma(\xi_f)^{-n(n-1)/2}\,p^{\epsilon}(\pi_f) }\right) 
=
\left(\frac
{p^{\epsilon^{\sigma}\cdot\epsilon_{\xi^{\sigma}}}(\pi_f^{\sigma}\otimes\xi_f^{\sigma})}
{\gamma(\xi_f^{\sigma})^{-n(n-1)/2}\,p^{\epsilon^{\sigma}}(\pi_f^{\sigma})}\right).
$$
\item Let ${\mathbb Q}(\pi_f, \xi_f)$ denote the compositum of the (number) fields
${\mathbb Q}(\pi_f)$ and ${\mathbb Q}(\xi_f)$. We have
$$
p^{\epsilon \cdot \epsilon_{\xi}}(\pi_f\otimes\xi_f) \,
\sim_{{\mathbb Q}(\pi_f, \xi_f)} \,
\gamma(\xi_f)^{-n(n-1)/2}\,p^{\epsilon}(\pi_f).
$$
By $\sim_{{\mathbb Q}(\pi_f, \xi_f)}$ we mean up to an element of ${\mathbb Q}(\pi_f, \xi_f)$.
\end{enumerate}
\end{thm}

\begin{proof}
Note that $(1) \Rightarrow (2)$ follows 
from the definitions of the rationality field of 
$\pi_f$ and $\xi_f$. It is really statement (1) which takes some work to prove; this
entails an analysis of the following diagram of maps. 
{\it Warning:} This diagram is not commutative! Indeed, the various
complex numbers involved in (1) measure the failure
of commutativity of this diagram. 

{\SMALL
\begin{equation}
\label{eqn:cube}
\xymatrix{
 & W(\pi_f)\ar[rr]^-{\mathcal{F}_{\pi_f , \epsilon}}\ar[d]^{W_{\xi_f}}\ar[ddl]_{\sigma} & &   
H^b(V_{\pi}\otimes M_{\mu}^{\vee})(\epsilon)
\ar[d]^{(A_{\xi}\otimes 1_{M_{\mu}^{\vee}})^*}\ar[ddl]_{\sigma} \\
 & W(\pi_f \otimes \xi_f) \ar[rr]^-{\mathcal{F}_{\pi_f\otimes\xi_f, \epsilon \cdot 
\epsilon_{\xi}}} \ar[ddl]_{\sigma} & & 
H^b(V_{\pi\otimes \xi} \otimes (M_{\mu}^{\vee}\otimes \xi_{\infty}^{-1}))
(\epsilon \cdot\epsilon_{\xi}) \ar[ddl]_{\sigma} \\
W(\pi^{\sigma}_f) \ar[rr]^-{\mathcal{F}_{\pi_f^{\sigma},\epsilon^{\sigma}}} 
\ar[d]_{W_{\xi_f^{\sigma}}} & & 
H^b(V_{\pi^{\sigma}} \otimes M_{\mu^{\sigma}}^{\vee})(\epsilon^{\sigma}) 
\ar[d]_{(A_{\xi^{\sigma}} \otimes 1_{M_{\mu^{\sigma}}^{\vee}})^*} & \\
W(\pi^{\sigma}_f \otimes \xi^{\sigma}_f) 
\ar[rr]^-{\mathcal{F}_{\pi_f^{\sigma} \otimes \xi_f^{\sigma}, 
\epsilon^{\sigma} \cdot\epsilon_{\xi^{\sigma}}}} & & 
H^b(V_{\pi^{\sigma} \otimes \xi^{\sigma}} \otimes 
(M_{\mu^{\sigma}}^{\vee}\otimes\xi_{\infty}^{-\sigma})) 
(\epsilon^{\sigma} \cdot\epsilon_{\xi^{\sigma}}) & 
}
\end{equation}}

We need to explain the undefined and abbreviated notations in the above diagram. 
We have abbreviated 
$H^b(\mathfrak{g}_{\infty}, K_{\infty}^0; V_{\pi}\otimes M_{\mu}^{\vee})(\epsilon)$ as
$H^b(V_{\pi}\otimes M_{\mu}^{\vee})(\epsilon)$. 
Same remark applies to three other objects. 
The maps $W_{\xi}$ and $A_{\xi}$ are defined as follows. If $w$ is any 
Whittaker function for ${\rm GL}_n({\mathbb A})$, then define 
$$
W_{\xi}(w)(g) = \xi({\rm det}(g))w(g)
$$ 
for all $g \in {\rm GL}_n({\mathbb A})$. It is easy to see that $W_{\xi}$ maps 
$W(\pi)$ onto $W(\pi \otimes \xi)$.
An identical formula defines $W_{\xi_f}$ and $W_{\xi_{\infty}}$.
Similarly, we define $A_{\xi}(\phi)$ for any automorphic form $\phi$ 
on ${\rm GL}_n({\mathbb A})$ by 
$$
A_{\xi}(\phi)(g) = \xi({\rm det}(g))\phi(g)
$$ for all 
$g \in {\rm GL}_n({\mathbb A})$.
It is easy
to see that $A_{\xi}$ maps $V_{\pi}$ onto $V_{\pi \otimes \xi}$.
The identity map 
on the vector space $M_{\mu}^{\vee}$ is denoted $1_{M_{\mu}^{\vee}}$. 
Observe that $A_{\xi} \otimes 1_{M_{\mu}^{\vee}}$ is a 
$G_{\infty}$-equivariant isomorphism from $V_{\pi}\otimes M_{\mu}^{\vee}$ onto
$V_{\pi\otimes \xi} \otimes (M_{\mu}^{\vee}\otimes \xi_{\infty}^{-1})$, 
and we denote $(A_{\xi} \otimes 1_{M_{\mu}^{\vee}})^*$ the map induced by 
$A_{\xi} \otimes 1_{M_{\mu}^{\vee}}$ in 
cohomology.

Analyzing the diagram involves an analysis of certain subdiagrams. Some of
these are independently interesting, and so we delineate them in the following
propositions. There are three propositions below which give information about
(non-)commutativity of some of the faces of (\ref{eqn:cube}).

\begin{prop}
\label{prop:whittaker-rational}
Let $\pi$ be a cuspidal automorphic representation of ${\rm GL}_n({\mathbb A})$. 
Let $\xi$ be a Hecke character of $F$. 
Then for any $\sigma \in {\rm Aut}({\mathbb C})$ we have
$$
\sigma \circ W_{\xi_f} = \sigma(\xi_f(t_{\sigma}^{-n(n-1)/2}))\,
W_{\xi_f^{\sigma}} \circ \sigma.
$$
\end{prop}

\begin{proof}
Consider the diagram
$$
\xymatrix{
W(\pi_f) \ar[d]^{W_{\xi_f}} \ar[r]^{\sigma}  &
W(\pi^{\sigma}_f) \ar[d]^{W_{\xi_f^{\sigma}}}  \\
W(\pi_f\otimes\xi_f) \ar[r]^{\sigma} & 
W(\pi^{\sigma}_f\otimes \xi^{\sigma}_f)
}
$$
Let $w \in W(\pi_f)$. We will chase this element $w$ in the above diagram. For all 
$g \in G({\mathbb A}_f)$ we have
\begin{eqnarray*}
\sigma(W_{\xi_f}(w))(g) & = & 
\sigma(W_{\xi_f}(w)(t_{\sigma, n} g)) \\
& = & \sigma(\xi_f({\rm det}(t_{\sigma, n} g)) w(t_{\sigma, n}g))  \\
& = & \sigma(\xi_f(t_{\sigma}^{-n(n-1)/2}))\sigma(\xi_f({\rm det}(g)))
\sigma(w(t_{\sigma, n}g)).
\end{eqnarray*}
On the other hand, we have 
\begin{eqnarray*}
W_{\xi_f^{\sigma}}(\sigma(w))(g) & = & 
\xi_f^{\sigma}({\rm det}(g)) \sigma(w)(g) \\
& = & \sigma(\xi_f({\rm det}(g)))\sigma(w(t_{\sigma,n}g)).
\end{eqnarray*}
Hence
$$
\sigma((W_{\xi_f}(w)) = \sigma(\xi_f(t_{\sigma}^{-n(n-1)/2}))
W_{\xi_f^{\sigma}}(\sigma(w)).
$$
(This measures the failure of commutativity of the above diagram.)
\end{proof}

\begin{lemma}
\label{lem:sigma-gauss}
$$
\sigma(\xi_f(t_{\sigma}^{-1})) = \frac{\sigma(\gamma(\xi_f))}{\gamma(\xi_f^{\sigma})}.
$$
\end{lemma}
\begin{proof}
Using the definition of $t_{\sigma}$ we have
\begin{eqnarray*}
\sigma(\gamma(\xi_f)) & = &\sigma(\gamma(\xi_f,y)) = 
\sum_{x \in (\mathcal{O}_F/\mathfrak{c})^{\times}} \sigma(\xi_f(x))
\sigma(e^{2 \pi i \, T_{F/{\mathbb Q}}(xy)}) \\
& = &
\sum_{x \in (\mathcal{O}_F/\mathfrak{c})^{\times}} \sigma(\xi_f(x))
e^{2 \pi i t_{\sigma}\, T_{F/{\mathbb Q}}(xy)}, 
\\
& = &  
\sum_{x \in (\mathcal{O}_F/\mathfrak{c})^{\times}} \sigma(\xi_f(x))
e^{2 \pi i \, T_{F/{\mathbb Q}}(t_{\sigma}xy)},   
\\
& = &  
\sum_{x \in (\mathcal{O}_F/\mathfrak{c})^{\times}} \sigma(\xi_f(t_{\sigma}^{-1}x))
e^{2 \pi i \, T_{F/{\mathbb Q}}(xy)} = 
\sigma(\xi_f(t_{\sigma}^{-1}))\gamma(\xi_f^{\sigma}).
\end{eqnarray*}
\end{proof}

\begin{cor}
\label{cor:whittaker-gauss}
$$
\sigma \circ W_{\xi_f} = 
\left(\frac{\sigma(\gamma(\xi_f))}{\gamma(\xi_f^{\sigma})}\right)^{n(n-1)/2} \,
W_{\xi_f^{\sigma}} \circ \sigma.
$$
\end{cor}

\begin{proof}
Follows from Proposition~\ref{prop:whittaker-rational} and Lemma~\ref{lem:sigma-gauss}. 
\end{proof}

\begin{prop}
\label{prop:automorphic-rational}
Let $\pi$ be a regular algebraic cuspidal automorphic representation of ${\rm 
GL}_n({\mathbb A})$, and let $\mu \in X^+(T_{\infty})$ 
be such that $\pi \in {\rm Coh}(G, \mu^{\vee})$. For any algebraic Hecke character 
$\xi$ we have 
$$
\sigma \circ (A_{\xi} \otimes 1_{M_{\mu}^{\vee}})^* = 
(A_{\xi^{\sigma}} \otimes 1_{M_{\mu^{\sigma}}^{\vee}})^* \circ \sigma.
$$
\end{prop}

\begin{proof}
This proof is a little more involved, and to help the reader, we adumbrate
it as follows. First go up to 
a bigger ambient vector space ($H^*_{dR}(\tilde{S},\mathcal{M})$)
and then use an identification of this bigger
space with another space ($H^*_{B}(\tilde{S},\mathcal{M})$) 
where it will be obvious that a lift of $A_{\xi}^*$ 
is Galois equivariant, and hence so is the original $A_{\xi}^*$. 
During the course of the proof, it helps to keep the following 
scheme of spaces and maps in mind. 
$$
\begin{array}{clcll}
H^b(\mathfrak{g}_{\infty},K_{\infty}^0; V_{\pi}\otimes M_{\mu}^{\vee})(\epsilon) 
& \simeq & H^b_{\rm cusp}(\tilde{S}, \mathcal{M}_{\mu}^{\vee})(\tilde{\pi})
& & \\
\downarrow
& & \downarrow & & \\
H^b(\mathfrak{g}_{\infty} ,K_{\infty}^0; 
\mathcal{A}_{\rm cusp}(G(F) \backslash G({\mathbb A})) 
\otimes M_{\mu}^{\vee})
& \simeq & H^b_{\rm cusp}(\tilde{S}, \mathcal{M}_{\mu}^{\vee})
& & \\
\downarrow & & \downarrow & & \\
H^b(\mathfrak{g}_{\infty},K_{\infty}^0; C^{\infty}(G(F) \backslash G({\mathbb A})) 
\otimes M_{\mu}^{\vee})
 & \simeq &
H^b_{dR}(\tilde{S}, \mathcal{M}_{\mu}^{\vee}) 
 & \simeq &
H^b_{B}(\tilde{S}, \mathcal{M}_{\mu}^{\vee}) 
\end{array}
$$
where all the vertical arrows are inclusions. Indeed, the rational structures on all
the above spaces comes from a rational structure on the Betti cohomology space on which
it is very easy to describe an action of ${\rm Aut}({\mathbb C})$ 
(see \cite[p.128]{clozel2}). In the above scheme, we need not (and so did not) refer to
cohomology with compact supports because one has \cite[p.129]{clozel2}
$$
H^*_{\rm cusp} \hookrightarrow H^*_c \to H^*_! := {\rm Image}(H^*_c \to H^*_{\rm dR}) 
$$
and the composite is an injection, i.e., $H^*_{\rm cusp} \hookrightarrow H^*_!$, and
hence $H^*_{\rm cusp} \hookrightarrow H^*_{\rm dR}$.

To begin the proof of Proposition~\ref{prop:automorphic-rational}, observe that
the map $(A_{\xi} \otimes 1_{M_{\mu}^{\vee}})^*$ 
in the statement of the proposition is the restriction 
to $H^b(\mathfrak{g}_{\infty},K_{\infty}^0; V_{\pi} \otimes 
M_{\mu}^{\vee})(\epsilon)$ of the map 
{\SMALL
$$
({}_{\mathcal{A}}A_{\xi} \otimes 1_{M_{\mu}^{\vee}})^* : 
H^b(\mathfrak{g}_{\infty},K_{\infty}^0; 
\mathcal{A}_{\rm cusp}(G(F) \backslash G({\mathbb A})) \otimes M_{\mu}^{\vee})  \to 
H^b(\mathfrak{g}_{\infty},K_{\infty}^0; 
\mathcal{A}_{\rm cusp}(G(F) \backslash G({\mathbb A})) 
\otimes (M_{\mu}^{\vee} \otimes \xi_{\infty}^{-1}))
$$}
induced by $A_{\xi}$ on $\mathcal{A}_{\rm cusp}(G(F) \backslash G({\mathbb A}))$. 
From a fundamental 
theorem of Borel \cite{borel}, cohomology with coefficients in the space of cusp forms
injects into cohomology with coefficients in the space of smooth functions, and the 
above map $({}_{\mathcal{A}}A_{\xi} \otimes 1_{M_{\mu}^{\vee}})^*$ is the restriction to 
$H^b(\mathfrak{g}_{\infty},K_{\infty}^0; 
\mathcal{A}_{\rm cusp}(G(F)\backslash G({\mathbb A})) \otimes M_{\mu}^{\vee})$ of the map
{\SMALL
$$
({}_{C^{\infty}}A_{\xi} \otimes 1_{M_{\mu}^{\vee}})^* : 
H^b(\mathfrak{g}_{\infty},K_{\infty}^0; 
C^{\infty}(G(F) \backslash G({\mathbb A})) \otimes M_{\mu}^{\vee}) \to
H^b(\mathfrak{g}_{\infty},K_{\infty}^0; 
C^{\infty}(G(F) \backslash G({\mathbb A})) \otimes 
(M_{\mu}^{\vee} \otimes \xi_{\infty}^{-1}))
$$}
induced by $A_{\xi}$ on $C^{\infty}(G(F) \backslash G({\mathbb A}))$.

On the other hand, relative Lie algebra cohomology can be identified with 
de Rham cohomology, and we can transport the map $({}_{C^{\infty}}A_{\xi} \otimes 
1_{M_{\mu}^{\vee}})^*$ across to get
$$
{}_{dR}A_{\xi}^* : H^b_{dR}(\tilde{S},\mathcal{M}_{\mu}^{\vee}) \to
H^b_{dR}(\tilde{S}, \mathcal{M}_{\mu}^{\vee} \otimes \xi_{\infty}^{-1}).
$$
(By $\mathcal{M}_{\mu}^{\vee} \otimes \xi_{\infty}^{-1}$ we mean the locally constant sheaf
induced by the representation $M_{\mu}^{\vee} \otimes \xi_{\infty}^{-1}$.)
We can describe the map ${}_{dR}A_{\xi}^*$ as follows. Let $K_f$ be an open compact subgroup
of $G({\mathbb A}_f)$ such that $\xi({\rm det}(k)) = 1$ for all $k \in K_f$.
Recall the identification (\cite[\S1.1]{labesse-schwermer})
$$
H^b_{dR}(S(K_f),\mathcal{M}_{\mu}^{\vee}) \simeq 
H^b(\mathfrak{g}_{\infty},K_{\infty}^0; C^{\infty}(G(F) \backslash G({\mathbb 
A}))^{K_f} \otimes M_{\mu}^{\vee}).
$$
The choice of $K_f$ implies that $A_{\xi}$ stabilizes 
$C^{\infty}(G({\mathbb Q}) \backslash G({\mathbb A}))^{K_f}$ and so induces a map 
$({}_{C^{\infty}}A_{\xi,K_f} \otimes 1_{M_{\mu}^{\vee}})^*$ 
on the right hand side above. Clearly, 
$\underrightarrow{\lim}_{K_f} \ 
({}_{C^{\infty}}A_{\xi,K_f} \otimes 1_{M_{\mu}^{\vee}})^* = 
({}_{C^{\infty}}A_{\xi} \otimes 1_{M_{\mu}^{\vee}})^*$.
Moving across to de Rham cohomology,
we can describe the map ${}_{dR}A_{\xi,K_f}^*$ as acting on differential forms
by {\it pointwise multiplication by $\xi$}, i.e., if $\omega$ is a 
$\mathcal{M}_{\mu}^{\vee}$-valued
(closed) differential form of degree $b$ on $S(K_f)$ then
$$
{}_{dR}A_{\xi,K_f}^*(\omega)_{\underline{g}} = \xi({\rm det}(g))\omega_{\underline{g}}
$$
for any $g \in {\rm GL}_n({\mathbb A})$, where $\underline{g}$ is the image of $g$ in 
$S(K_f)$. (For any $x \in S(K_f)$, 
by $\omega_x$ we mean the value at $x$ of $\omega$ which is a section of the 
$b$-th exterior of the cotangent bundle twisted by $\mathcal{M}_{\mu}^{\vee}$ over 
the manifold $S(K_f)$.) Observe that the above equation is
well defined. Passing to the limit we get the map
${}_{dR}A_{\xi}^* = \underrightarrow{\rm lim} \  {}_{dR}A_{\xi,K_f}^*.$

Now we move across to Betti cohomology via the de Rham isomorphism, and get a map 
$$
{}_BA_{\xi}^* : H^b_B(\tilde{S}, \mathcal{M}_{\mu}^{\vee}) \to
H^b_B(\tilde{S}, \mathcal{M}_{\mu}^{\vee}\otimes\xi_{\infty}^{-1}).
$$
The point of going to Betti cohomology is because the action of ${\rm Aut}({\mathbb C})$
is especially simple to describe--it acts by acting on the coefficient system. 
(See \cite[page 128]{clozel2}.) Namely, if $\sigma \in {\rm Aut}({\mathbb C})$
then we have a $\sigma$ linear isomorphism
$$
H^*_B(S(K_f),\mathcal{M}_{\mu}^{\vee}) \to H^*_B(S(K_f), 
\mathcal{M}_{\mu^{\sigma}}^{\vee}).
$$
This isomorphism is the one induced in cohomology by the following map on the singular cochain
complex
$$
\Delta^*(S(K_f), \mathcal{M}_{\mu}^{\vee}) \to 
\Delta^*(S(K_f), \mathcal{M}_{\mu^{\sigma}}^{\vee})
$$
which is simply $\tau \mapsto l_{\sigma} \circ \tau$, if $l_{\sigma}$ is the 
$\sigma$-linear isomorphism from $M_{\mu}^{\vee}$ to $M_{\mu^{\sigma}}^{\vee}$.
(Recall that $M_{\mu}^{\vee}$ is 
defined over ${\mathbb Q}$ and that the action of $\sigma$ on $\mu$ is via the 
embeddings of $F$ into ${\mathbb C}$.)
The action of ${\rm Aut}({\mathbb C})$ on 
$H^*_B(S(K_f),\mathcal{M}_{\mu}^{\vee})$ 
can be transported to an action on 
$H^*_{dR}(S(K_f),\mathcal{M}_{\mu}^{\vee})$, and after passing to the limit, 
induces an action of ${\rm Aut}({\mathbb C})$ on each of the spaces
$$
H^*(\mathfrak{g}_{\infty},K_{\infty}^0, V_{\pi}\otimes M_{\mu}^{\vee})(\epsilon) 
\subset
H^*_{\rm cusp}(\tilde{S},\mathcal{M}_{\mu}^{\vee}) \subset 
H^*_{dR}(\tilde{S},\mathcal{M}_{\mu}^{\vee}).
$$

The statement in the proposition may be phrased as that the map 
$(A_{\xi} \otimes 1_{M_{\mu}^{\vee}})^*$ is 
${\rm Aut}({\mathbb C})$-equivariant. From the above description of the 
action of ${\rm Aut}({\mathbb C})$, we can see that 
the $(A_{\xi} \otimes 1_{M_{\mu}^{\vee}})^*$ is ${\rm Aut}({\mathbb C})$-equivariant
if and only if 
any of the maps $({}_{\mathcal{A}}A_{\xi} \otimes 1_{M_{\mu}^{\vee}})^*$,
$({}_{C^{\infty}}A_{\xi} \otimes 1_{M_{\mu}^{\vee}})^*$,  ${}_{dR}A_{\xi}^*$ or 
${}_{B}A_{\xi}^*$ is 
${\rm Aut}({\mathbb C})$-equivariant. 

It is easy to see that the map ${}_{dR}A_{\xi}^*$ is ${\rm Aut}({\mathbb C})$-equivariant,
since ${\rm Aut}({\mathbb C})$ acts on the coefficient system, 
and ${}_{dR}A_{\xi}^*$ is more intrinsic to the manifold. 
More precisely, consider the de Rham map 
$\Omega^*(S(K_f), \mathcal{M}^{\vee}_{\mu}) \to
\Delta^*(S(K_f), \mathcal{M}^{\vee}_{\mu})$ from the 
space of differential forms 
to the space of singular cochains, given by integration. 
(This induces the de Rham isomorphism in cohomology.)
We can describe the action of $\sigma \in {\rm Aut}({\mathbb C})$ on 
$\omega \in \Omega^b(S(K_f), \mathcal{M}^{\vee}_{\mu})$,
by $\sigma(\omega)_x = l_{\sigma} \circ \omega_x$ for $x \in S(K_f)$. 
For any $g \in {\rm GL}_n({\mathbb A})$, if $\underline{g}$ denotes the image of $g$  
in $S(K_f)$, we have
\begin{eqnarray*}
\sigma({}_{dR}A_{\xi, K_f}(\omega))_{\underline{g}}
& = &
  l_{\sigma} \circ  {}_{dR}A_{\xi, K_f}(\omega)_{\underline{g}} = 
  l_{\sigma} \circ \xi({\rm det}(g))\omega_{\underline{g}} \\
& = &  \sigma(\xi({\rm det}(g))) l_{\sigma} \circ \omega_{\underline{g}} = 
{}_{dR}A_{\xi^{\sigma}, K_f}(\sigma(\omega))_{\underline{g}}.
\end{eqnarray*}
In other words, $\sigma \circ {}_{dR}A_{\xi, K_f} = 
{}_{dR}A_{\xi^{\sigma}, K_f} \circ \sigma$. Passing to the limit over all $K_f$, we
get $\sigma \circ {}_{dR}A_{\xi} = 
{}_{dR}A_{\xi^{\sigma}} \circ \sigma$, which induces the desired equality of maps
in cohomology. 
\end{proof}

\begin{prop}
\label{prop:diagram}
The diagram
$$
\xymatrix{
W(\pi_f) \ar[d]^{W_{\xi_f}} \ar[rrr]^-{\mathcal{F}_{\pi_f, \epsilon}} & & &
H^b(\mathfrak{g}_{\infty},K_{\infty}^0; V_{\pi} \otimes M_{\mu}^{\vee})(\epsilon) 
\ar[d]^{(A_{\xi} \otimes 1_{M_{\mu}^{\vee}})^*} \\
W(\pi_f\otimes\xi_f) \ar[rrr]^-{\mathcal{F}_{\pi_f\otimes \xi_f, 
\epsilon \cdot \epsilon_{\xi}}} & & &
H^b(\mathfrak{g}_{\infty},K_{\infty}^0; V_{\pi \otimes \xi} \otimes 
(M_{\mu}^{\vee}\otimes\xi_{\infty}^{-1}))(\epsilon \cdot \epsilon_{\xi})
}
$$
commutes. (The horizontal maps are the {\it un-normalized} maps.)
\end{prop}

\begin{proof}
That this diagram commutes can be seen by observing that the following three diagrams 
commute, since the horizontal maps are both 
the compositum of three maps. (For brevity we denote $H^b(-)$ 
for $H^b(\mathfrak{g}_{\infty},K_{\infty}^0; -)$.)
{\Small
\begin{equation}
\label{eqn:diagram1}
\xymatrix{
W(\pi_f) \ar[d]^{W_{\xi_f}} \ar[r] &
W(\pi_f) \otimes 
H^b(W(\pi_{\infty})\otimes M_{\mu}^{\vee})(\epsilon) 
\ar[d]^{W_{\xi_f}\otimes (W_{\xi_{\infty}}\otimes1_{M_{\mu}^{\vee}})^*} \\
W(\pi_f\otimes\xi_f) \ar[r] &  
W(\pi_f \otimes \xi_f) \otimes 
H^b(W(\pi_{\infty}\otimes\xi_{\infty}) 
\otimes (M_{\mu}^{\vee} \otimes\xi_{\infty}^{-1}))(\epsilon\cdot\epsilon_{\xi})
}
\end{equation}}

{\Small
\begin{equation}
\label{eqn:diagram2}
\xymatrix{
W(\pi_f) \otimes H^b(W(\pi_{\infty})\otimes M_{\mu}^{\vee})(\epsilon)
\ar[d]^{W_{\xi_f}\otimes (W_{\xi_{\infty}}\otimes1_{M_{\mu}^{\vee}})^*} 
\ar[r] &
H^b(W(\pi) \otimes M_{\mu}^{\vee})(\epsilon) 
\ar[d]^{(W_{\xi}\otimes1_{M_{\mu}^{\vee}})^*}\\
W(\pi_f\otimes\xi_f) \otimes
H^b(W(\pi_{\infty}\otimes\xi_{\infty})\otimes (M_{\mu}^{\vee}\otimes\xi_{\infty}^{-1}))
(\epsilon\cdot\epsilon_{\xi}) \ar[r] &
H^b(W(\pi \otimes \xi) \otimes (M_{\mu}^{\vee}\otimes\xi_{\infty}^{-1}))
(\epsilon\cdot\epsilon_{\xi})
}
\end{equation}
}

{\Small
\begin{equation}
\label{eqn:diagram3}
\xymatrix{
H^b(W(\pi)\otimes M_{\mu}^{\vee})(\epsilon) 
\ar[d]^{(W_{\xi}\otimes1_{M_{\mu}^{\vee}})^*} 
\ar[r] &
H^b(V_{\pi} \otimes M_{\mu}^{\vee})(\epsilon)
\ar[d]^{(A_{\xi}\otimes1_{M_{\mu}^{\vee}})^*} \\
H^b(W(\pi \otimes \xi)\otimes (M_{\mu}^{\vee}\otimes\xi^{-1}))(\epsilon\cdot\epsilon_{\xi}) 
\ar[r] &
H^b(V_{\pi \otimes \xi} \otimes 
(M_{\mu}^{\vee}\otimes\xi_{\infty}^{-1}))(\epsilon\cdot\epsilon_{\xi})
}
\end{equation}
}
For the commutativity of (\ref{eqn:diagram1}), note that the linear map 
$W_{\xi_{\infty}}$ induces a $G_{\infty}$-equivariant isomorphism 
$W_{\xi_{\infty}}\otimes 1_{M_{\mu}^{\vee}} : W(\pi_{\infty})\otimes M_{\mu}^{\vee} \to
W(\pi_{\infty}\otimes\xi_{\infty})\otimes (M_{\mu}^{\vee} \otimes \xi_{\infty}^{-1})$, and 
hence induces an isomorphism $(W_{\xi_{\infty}}\otimes 1_{M_{\mu}^{\vee}})^*$ in cohomology.
From the choice we made on the generators of such one dimensional cohomology spaces we have 
$(W_{\xi_{\infty}}\otimes 1_{M_{\mu}^{\vee}})^*({\bf w}(\pi_{\infty},\epsilon)) = 
{\bf w}(\pi_{\infty}\otimes\xi_{\infty},\epsilon\cdot \epsilon_{\xi})$. Now it is easy to 
check that (\ref{eqn:diagram1}) commutes. 
The diagram in 
(\ref{eqn:diagram2}) is the one induced in cohomology of the commuative diagram 
$$
\xymatrix{
W(\pi_f) \otimes W(\pi_{\infty}) \otimes M_{\mu}^{\vee}
\ar[r] \ar[d]^{W_{\xi_f}\otimes W_{\xi_{\infty}} \otimes 1_{M_{\mu}^{\vee}}}
 & W(\pi) \otimes M_{\mu}^{\vee}
\ar[d]^{W_{\xi}\otimes 1_{M_{\mu}^{\vee}}} \\
W(\pi_f\otimes\xi_f) \otimes W(\pi_{\infty} \otimes\xi_{\infty}) 
\otimes (M_{\mu}^{\vee} \otimes \xi_{\infty}^{-1})
\ar[r] & W(\pi\otimes\xi) \otimes (M_{\mu}^{\vee} \otimes \xi_{\infty}^{-1})
}
$$
and hence (\ref{eqn:diagram2}) commutes. Finally, the diagram in 
(\ref{eqn:diagram3}) is the one induced in cohomology of the diagram 
$$
\xymatrix{
W(\pi) \otimes M_{\mu}^{\vee} \ar[r] 
\ar[d]^{W_{\xi} \otimes 1_{M_{\mu}^{\vee}}} & 
V_{\pi} \otimes M_{\mu}^{\vee} \ar[d]^{A_{\xi} \otimes 1_{M_{\mu}^{\vee}}} \\
W(\pi\otimes\xi) \otimes (M_{\mu}^{\vee} \otimes \xi_{\infty}^{-1})
\ar[r] & 
V_{\pi\otimes\xi} \otimes (M_{\mu}^{\vee} \otimes \xi_{\infty}^{-1})
}
$$
where the horizontal maps are the inverses of taking the Fourier coefficients. If the 
Fourier coeffcient map $V_{\pi} \to W(\pi) = W(\pi,\psi)$ is denoted 
$\phi \mapsto w_{\psi}(\phi)$, where
$$
w_{\psi}(\phi)(g) = 
\int_{N(F) \backslash N({\mathbb A})}
\phi(ng) \overline{\psi(n)} \, dn
$$
then it is easy to see that $w_{\psi}(\phi\otimes\xi) = w_{\psi}(\phi)\otimes\xi$, hence
this last diagram commutes.
\end{proof}

We can now finish the proof of Theorem~\ref{thm:twisted} as follows. We consider the 
composite map $(A_{\xi^{\sigma}} \otimes 1_{M_{\mu^{\sigma}}^{\vee}})^* \circ \sigma \circ 
\mathcal{F}_{\pi_f,\epsilon}$ 
in the diagram (\ref{eqn:cube}). On the one hand, using 
Equation~\ref{eqn:unnormalized} and Proposition~\ref{prop:diagram}, we have
\begin{eqnarray*}
(A_{\xi^{\sigma}} \otimes 1_{M_{\mu^{\sigma}}^{\vee}})^*
\circ \sigma \circ \mathcal{F}_{\pi_f,\epsilon}
& = & (A_{\xi^{\sigma}} \otimes 1_{M_{\mu^{\sigma}}^{\vee}})^* \circ 
      \left(\frac{\sigma(p^{\epsilon}(\pi_f)}{p^{\epsilon^{\sigma}}(\pi_f^{\sigma})}\right)
      \mathcal{F}_{\pi_f^{\sigma},\epsilon^{\sigma}} \circ \sigma \\
& = & \left(\frac{\sigma(p^{\epsilon}(\pi_f))}{p^{\epsilon^{\sigma}}(\pi_f^{\sigma})}\right)
      \mathcal{F}_{\pi_f^{\sigma}\otimes\xi_f^{\sigma},\epsilon^{\sigma}\cdot\epsilon_{\xi^{\sigma}}}
\circ W_{\xi_f^{\sigma}} \circ \sigma.
\end{eqnarray*}
On the other hand, using Propositions~\ref{prop:automorphic-rational},
\ref{prop:diagram}, Equation~\ref{eqn:unnormalized}
and Corollary~\ref{cor:whittaker-gauss} (in that order), we have 
{\Small
\begin{eqnarray*}
(A_{\xi^{\sigma}} \otimes 1_{M_{\mu^{\sigma}}^{\vee}})^*
\circ \sigma \circ \mathcal{F}_{\pi_f,\epsilon}
& = & 
\sigma \circ (A_{\xi} \otimes 1_{M_{\mu}^{\vee}})^* \circ \mathcal{F}_{\pi_f,\epsilon}\\
& = & \sigma \circ \mathcal{F}_{\pi_f\otimes\xi_f, \epsilon\cdot\epsilon_{\xi}}
\circ W_{\xi_f} \\
& = & \left(\frac
 {\sigma(p^{\epsilon\cdot\epsilon_{\xi}}(\pi_f\otimes\xi_f))}
 {p^{ \epsilon^{\sigma} \cdot \epsilon_{\xi^{\sigma}} }
 (\pi_f^{\sigma} \otimes \xi_f^{\sigma})}\right)
 \mathcal{F}_{\pi_f^{\sigma}\otimes\xi_f^{\sigma},\epsilon^{\sigma}\cdot\epsilon_{\xi^{\sigma}}} 
\circ \sigma \circ W_{\xi_f} 
\\
& = & 
\left(\frac 
 {\sigma(p^{\epsilon\cdot\epsilon_{\xi}}(\pi_f\otimes\xi_f))}
 {p^{ \epsilon^{\sigma} \cdot \epsilon_{\xi^{\sigma}} }
 (\pi_f^{\sigma} \otimes \xi_f^{\sigma})}\right)
 \left(\frac{\sigma(\gamma(\xi_f))}{\gamma(\xi_f^{\sigma})}\right)^{n(n-1)/2}
 \mathcal{F}_{\pi_f^{\sigma}\otimes\xi_f^{\sigma},\epsilon^{\sigma}\cdot\epsilon_{\xi^{\sigma}}} 
 \circ W_{\xi_f^{\sigma}} \circ \sigma.
\end{eqnarray*}}
Putting both together we have
$$
\frac{\sigma(p^{\epsilon}(\pi_f))}{p^{\epsilon^{\sigma}}(\pi_f^{\sigma})}
= 
\left(\frac
 {\sigma(p^{\epsilon\cdot\epsilon_{\xi}}(\pi_f\otimes\xi_f))}
 {p^{ \epsilon^{\sigma} \cdot \epsilon_{\xi^{\sigma}} }
 (\pi_f^{\sigma} \otimes \xi_f^{\sigma})}\right)
 \left(\frac{\sigma(\gamma(\xi_f))}{\gamma(\xi_f^{\sigma})}\right)^{n(n-1)/2}
$$
from which the theorem follows. 
\end{proof}

\subsection{Some remarks}% and an example}
\label{sec:remarks}

\begin{rem}{\rm
The reader should compare Theorem~\ref{thm:twisted} with the conjectures of Blasius 
\cite{blasius2} 
and Panchiskin \cite{panchiskin}, on the behaviour of Deligne's motivically defined 
periods upon twisting a given motive by Artin motives. (A finite order Hecke character is an example
of an Artin motive.)}
\end{rem}

\begin{rem}{\rm
Suppose $\pi$ is an algebraic cuspidal representation of ${\rm GL}_n({\mathbb 
A}_F)$, and suppose $M = M(\pi)$ is the conjectural motive attached to $\pi$ as in 
\cite[Conjecture 
4.5]{clozel2}. A very interesting question (modulo this conjecture) is to understand how the 
periods
$p^{\epsilon}(\pi_f)$ considered above compare with Deligne's periods $c^{\pm}(M)$ attached to $M$.
(See \cite{deligne}.)
This is related to the issue of factorization of Deligne's periods over the infinite places as
in Blasius \cite[M.8]{blasius2}. The question is even more delicate than just a factorization issue,
because according to \cite[Conjecture 2.8]{deligne}, Deligne's periods are meant to capture the 
transcendental part of the special values of the
motivic $L$-function $L(s,M)$, however, Mahnkopf's work \cite{mahnkopf2} suggests that
the transcendental part of the special values of 
$L(s,\pi_f)$ depends not only on the periods $p^{\epsilon}(\pi_f)$, but also on the periods 
attached to 
a sequence of representations $\pi_i$ of ${\rm GL}_{n-2i}({\mathbb A}_F)$ for $1 \leq i < n/2$.}
\end{rem}

\begin{exam}[Special case of Tate's conjecture]
{\rm
Consider Theorem~\ref{thm:twisted} in the following special case: Let $F$ be a real quadratic 
extension of ${\mathbb Q}$. Let $K/F$ be a totally imaginary quadratic extension. Let $\chi$ be 
a Hecke character of $K$ which is ${\rm Gal}(K/F)$-regular, and let $\pi = {\rm AI}_{K/F}(\chi)$
be the cuspidal automorphic representation of ${\rm GL}_2({\mathbb A}_F)$ obtained by automorphic
induction. Assume that the infinity type of $\chi$ is so arranged that $\pi$ is algebraic and 
regular. Note that $\pi \simeq \pi\otimes \omega_{K/F}$, where $\omega_{K/F}$ is the quadratic 
Hecke character of $F$ associated to $K/F$ by class field theory. 
If we denote a signature $\epsilon = (\epsilon_1,\epsilon_2)$ 
with $\epsilon_i = \pm$ (instead of $1\!\!1$ or sgn as before), then 
the signature of $\omega_{K/F}$ is $(-,-)$. In this setting, Theorem~\ref{thm:twisted}
gives
$$
p^{++}(\pi_f)/p^{--}(\pi_f) \sim \gamma(\omega_{K/F}), \ \ 
p^{+-}(\pi_f)/p^{-+}(\pi_f) \sim \gamma(\omega_{K/F}).
$$
A Hilbert modular form over $F$ of CM type corresponds to such a representation and the above period
relations are formally the same relations proved by Murty and Ramakrishnan.
(See \cite[Theorem A]{kumar-dinakar}.) Such a period relation is one of the main 
ingredients for them to prove Tate's conjecture in a special case.
}
\end{exam}

\section{Deligne's conjecture}
\label{sec:dc}

\subsection{Statement of the conjecture}
\label{sec:deligne}

Let $\varphi \in S_k(N,\omega)$, i.e., $\varphi$ is a 
holomorphic cusp form on the upper half plane,
for $\Gamma_0(N)$, of weight $k$ and nebentypus character
$\omega$. Let $\varphi(z) = \sum_{n=1}^{\infty} a_nq^n$ be the Fourier 
expansion of $\varphi$ at infinity. 
We let $L(s, \varphi)$ stand for the completed $L$-function associated to
$\varphi$ and let $L_f(s, \varphi)$ stand for its finite part. 
Assume that $\varphi$ is a primitive form in $S_k(N,\omega)$.
By primitive, we mean that it is an eigenform, a newform and is normalized
such that $a_1(\varphi)=1$. 
For ${\rm Re}(s) \gg 0$, the finite part $L_f(s, \varphi)$ is 
a Dirichlet series with an Euler product as
$$
L_f(s, \varphi) = 
\sum_{n=1}^{\infty} a_n n^{-s} = 
\prod_p L_p(s, \varphi)
$$
where, for all primes $p$, we have 
$$
L_p(s, \varphi) = 
(1-a_pp^{-s} + \omega(p)p^{k-1-2s})^{-1} = 
(1-\alpha_{p,\varphi} p^{-s})^{-1} (1-\beta_{p,\varphi} p^{-s})^{-1}
$$
with the convention that if $p|N$ then $\beta_{p,\varphi} = 0$. We let 
${\rm Supp}(N)$ stand for the set of primes dividing $N$ and let 
$S = {\rm Supp}(N)\cup \{\infty\}$. 

For any $n \geq 1$, the partial $n$-th symmetric power $L$-function
is defined as 
$$
L^S(s, {\rm Sym}^n \varphi) = 
\prod_{p \notin S} L_p(s, {\rm Sym}^n \varphi ), \ \ {\rm Re}(s) \gg 0,
$$
where, for all $p \notin S$, we have 
$$
L_p(s,{\rm Sym}^n \varphi ) = 
\prod_{i=0}^n (1-\alpha_{p,\varphi}^i \beta_{p,\varphi}^{n-i} 
p^{-s})^{-1}.
$$

The Langlands program predicts that $L^S(s, {\rm Sym}^n \varphi)$, which 
is intially defined only in a half plane, admits a meromorphic 
continuation to the entire complex plane and that it has all the usual 
properties an automorphic $L$-function is supposed to have. This is known 
for $n \leq 9$ (with only partial results for $n \geq 5$)
from the works of several people including Hecke, Shimura, 
Jacquet--Langlands, Gelbart--Jacquet, Kim and Shahidi. 
It is also known for all $n$    
for cusp forms of a special 
type, for instance, if the representation corresponding to the cusp form 
is dihedral or the other polyhedral types. 

For $\varphi$ a primitive form in $S_k(N,\omega)$, we
let $M(\varphi)$ be the motive associated to $\varphi$. This is a rank
two motive over ${\mathbb Q}$ with coefficients in the field
${\mathbb Q}(\varphi)$ generated by the Fourier coefficients of $\varphi$.
(We refer the reader to Deligne \cite{deligne} and 
Scholl \cite{scholl} for details about $M(\varphi)$.)
The $L$-function $L(s, M(\varphi))$ associated to this motive is
$L(s, \varphi)$.
Given the motive $M(\varphi)$, there are nonzero complex numbers,
called Deligne's periods, $c^{\pm}(M(\varphi))$ associated to it.
Similarly, for the symmetric powers
${\rm Sym}^n(M(\varphi))$, we have the corresponding periods
$c^{\pm}({\rm Sym}^n(M(\varphi)))$. In \cite[Proposition 7.7]{deligne}
the periods for the symmetric powers are related to the
periods of $M(\varphi)$. The explicit formulae therein have
a quantity
$\delta(M(\varphi))$, which depends on the Gauss sum of the nebentypus
character $\omega$ and the weight $k$, and is given by
$$
\delta(M(\varphi)) \sim
(2 \pi i)^{1-k}\gamma(\omega) :=
(2 \pi i)^{1-k}
\sum_{u=0}^{c-1}
\omega(u) {\rm exp}(-2\pi i u/c),
$$
where $c$ is the conductor of $\omega$, and 
by $\sim$ we mean up to an element of ${\mathbb Q}(\varphi)$. 
We will denote the right hand side by $\delta(\omega)$. For brevity, we 
will denote $c^{\pm}({\rm Sym}^n(M(\varphi)))$ by $c^{\pm}({\rm Sym}^n(\varphi))$, 
and if $n=1$ this will be denoted $c^{\pm}(\varphi)$.

Recall \cite[Definition 1.3]{deligne} that an integer $m$ is 
{\it critical} for
any motivic $L$-function $L(s,M)$ if
both $L_{\infty}(s,M)$ and $L_{\infty}(1-s, M^{\vee})$ are regular at 
$s=m$. Recall also that given a $\varphi$ as above, and given any 
$\sigma \in {\rm Aut}({\mathbb C})$, the function
$\varphi^{\sigma}(z) = \sum_{n=1}^{\infty}\sigma(a_n)e^{2\pi inz}$ is  
a primitive modular form in $S_k(N, \omega^{\sigma})$.
With this notation, 
we now state Deligne's conjecture 
\cite[Section 7]{deligne} on the special values of the symmetric power
$L$-functions.

\begin{con}
\label{conj:deligne}
Let $\varphi$ be a primitive form in $S_k(N,\omega)$.
There exist nonzero complex numbers $c^{\pm}(\varphi)$ such that 
\begin{enumerate}
\item for $n = 2r+1$, if we denote $d^{\pm} = r+1$, and  
$$
c^{\pm}({\rm Sym}^n\varphi) = 
c^{\pm}(\varphi)^{(r+1)(r+2)/2}\,
c^{\mp}(\varphi)^{r(r+1)/2}\, 
\delta(\omega)^{r(r+1)/2}; 
$$
\item for $n = 2r$, if we denote $d^+ = r+1$, $d^- = r$, 
\begin{eqnarray*}
c^+({\rm Sym}^n\varphi) & = & (c^+(\varphi)c^-(\varphi))^{r(r+1)/2}\,
\delta(\omega)^{r(r+1)/2},\ \ {\rm and} \\
c^-({\rm Sym}^n\varphi) & = &(c^+(\varphi)c^-(\varphi))^{r(r+1)/2}\,
\delta(\omega)^{r(r-1)/2};
\end{eqnarray*}
\end{enumerate}
then for all $\sigma \in {\rm Aut}({\mathbb C})$ and for any integer $m$ which 
is critical for $L_f(s, {\rm Sym}^n \varphi)$ we have
$$
\sigma \left(
\frac{L_f(m, {\rm Sym}^n\varphi)}
{(2\pi i)^{m d^{\pm}}\, c^{\pm}({\rm Sym}^n \varphi)}\right) = 
\frac{L_f(m, {\rm Sym}^n \varphi^{\sigma})}
{(2\pi i)^{md^{\pm}}\, c^{\pm}({\rm Sym}^n \varphi^{\sigma})}\,
$$
where $\pm = (-1)^m$.

\end{con}

We wish to emphasize that in the original conjecture of Deligne, 
the numbers $c^{\pm}$ are {\it periods} which come via a comparison of rational  
structures on the de Rham and Betti realization of the motive; however,
in this paper they are just a couple of complex numbers in terms of which the 
critical values of the symmetric power $L$-functions can be expressed. 

For $n \leq 3$ the conjecture is known. For $n=1$ it is due to Shimura \cite{shimura2}; 
for $n=2$ it is due to Sturm \cite{sturm1}, \cite{sturm2}; for $n=3$ it is due
to Garrett--Harris \cite{garrett-harris} and Kim--Shahidi \cite{kim-shahidi-israel}. 
In general the conjecture is not known for higher ($n \geq 4$) symmetric
power $L$-functions. Although, if $\varphi$ is dihedral, then
the conjecture is known to be true for any symmetric power via 
motivic techniques. This is because Deligne's main conjecture 
\cite[\S1 and \S2]{deligne} is known 
if one considers only the motives as those attached to abelian varieties
and the category used is that defined by using absolute Hodge cycles for
morphisms. However, in \S\ref{sec:dihedral} we give a proof 
in the dihedral case using only $L$-functions. 

We remark that a prelude to this conjecture was certain calculations 
made by Zagier \cite{zagier} wherein he showed that such a statement 
holds for the $n$-th symmetric power $L$-function, with $n \leq 4$, of 
the Ramanujan $\Delta$-function.

\subsection{Some known results on special values of $L$-functions}
\label{sec:potpourri}

The purpose of this subsection is to recall some known special values 
results which will be relevant to us. Our aim in \S\ref{sec:dihedral} is to prove
a special values theorem for $L_f(m, {\rm Sym}^n \varphi)$ when $\varphi$ is 
dihedral, in which case this $L$-function factorizes into a product of 
$L$-functions of degrees at most two; see Corollary \ref{cor:dihedral}.
The following theorems describe 
the special values of these $L$-functions. In the course of the proof, we will need
to use the special values of symmetric squared $L$-functions, which is also 
recalled in this subsection.

\begin{thm}[Dirichlet $L$-functions]
\label{thm:dirichlet}
Let $\chi$ be a nontrivial, primitive Dirichlet character modulo $N$. 
Let $L_f(\chi, s) = \sum_{n=1}^{\infty} \chi(n)/n^s$ be the usual 
Dirichlet $L$-series associated to $\chi$.
Let $\nu \in \{0,1\}$ be defined by $\chi(-1)=(-1)^{\nu}$. For any integer 
$m \geq 1$, with $m \equiv \nu \pmod{2}$, we have 
$$
L_f(m, \chi) = (-1)^{1+(m-\nu)/2}
\frac{\gamma(\chi)}{2i^{\nu}}
\left(\frac{2\pi}{N}\right)^m \frac{B_{m,\overline{\chi}}}{m!}.
$$
\end{thm}

In the above theorem, a proof of which may be found in \cite[\S VII.9]{neukirch},
the quantity $B_{m,\chi}$ is a generalized Bernoulli number which lies 
in ${\mathbb Q}(\chi)$--the field generated by 
the values of $\chi$--and is ${\rm Aut}({\mathbb C})$-equivariant. 
For our purposes we rephrase the above theorem as 
\begin{equation}
\label{eqn:dirichlet}
\sigma\left(\frac{L_f(m, \chi)}{(2\pi i)^m \gamma(\chi)}\right) = \\
\frac{L_f(m, \chi^{\sigma})}{(2\pi i)^m \gamma(\chi^{\sigma})}.
\end{equation}
Even if $\chi$ is not primitive, the above equation still holds. Suppose that
$\chi$ is a character modulo $N$ and of conductor $c$, with $\chi_0$ the 
associated primitive character modulo $c$. Then, by definition, we have
$\gamma(\chi) = \gamma(\chi_0)$, and moreover, if
$S$ is the (possibly empty) set of primes $p$ which divides $N$ but
not $c$, then we have the equality
$$
L_f(s,\omega) =
L_f(s,\omega_0)\prod_{p \in S}(1-\omega_0(p)p^{-s})
$$
in a half plane and hence everywhere. It is now easy to see that 
(\ref{eqn:dirichlet}) holds for such a possibly imprimitive $\chi$.

\begin{thm}[Modular forms; Shimura \cite{shimura1}, \cite{shimura2}]
\label{thm:shimura}
Let $\varphi$ be a primitive cusp form in $S_k(N,\omega)$
with Fourier expansion $\sum_{n=1}^{\infty} a_nq^n$ where $q = e^{2\pi i 
z}$. Let ${\mathbb Q}(\varphi)$ be the 
field generated over ${\mathbb Q}$ by the Fourier coefficients of 
$\varphi$. Let $\xi$ be a Dirichlet character and 
let $L_f(s,\varphi,\xi) = \sum_{n=1}^{\infty}\xi(n)a_nn^{-s}.$
Let $m$ be an integer with $1 \leq m \leq k-1$.
There exists complex numbers $u^{\pm}(\varphi)$
depending only on $\varphi$ such that
\begin{enumerate}

%\item Let $\rho$ denote complex conjugation. 
%For all $\sigma \in {\rm Aut}({\mathbb C})$ we have
%$$
%u^{\pm}(\varphi^{\sigma})^{\rho} = \pm u^{\pm}(\varphi^{\sigma \rho}).
%$$

\item For all $\sigma \in {\rm Aut}({\mathbb C})$ we have
$$
\sigma\left(\frac{L_f(m,\varphi,\xi)}
{(2\pi i)^m \, \gamma(\xi) \, u^{\pm}(\varphi)} \right) = 
\frac{L_f(m,\varphi^{\sigma}, \xi^{\sigma})}
{(2\pi i)^m \, \gamma(\xi^{\sigma}) \, u^{\pm}(\varphi^{\sigma})}
$$
where $\pm = (-1)^m\xi(-1)$.

\item Let $\langle \varphi , \varphi \rangle$ 
be the Petersson inner product defined as in 
\cite[(2.1)]{shimura1}. For all $\sigma \in {\rm Aut}({\mathbb C})$ we have
$$
\sigma \left(\frac{i^{1-k}\pi \gamma(\omega)
\langle \varphi , \varphi \rangle} {u^+(\varphi)u^-(\varphi)}\right) \ =\  
\frac{i^{1-k}\pi \gamma(\omega^{\sigma})
\langle \varphi^{\sigma} , \varphi^{\sigma} \rangle} {u^+(\varphi^{\sigma})u^-(\varphi^{\sigma})}.
$$
\end{enumerate}
\end{thm}

Some remarks are in order, especially about the Shimura's periods 
$u^{\pm}(\varphi)$ and their relation to Deligne's periods 
$c^{\pm}(\varphi)$ of \S\ref{sec:deligne}. 
If $k \geq 3$, then Shimura's periods are defined as 
$$ 
u^+(\varphi) = \frac{L_f(k-1, \varphi, \xi^+)}
{(2\pi i)^{k-1} \, \gamma(\xi^+)}, \ \  
u^-(\varphi) = \frac{L_f(k-1, \varphi, \xi^-)}
{(2\pi i)^{k-1} \, \gamma(\xi^-)}
$$
where $\xi^{\pm}$ are fixed real valued characters such that
$\xi^+(-1) = (-1)^{k-1}$ and $\xi^-(-1) = (-1)^k$.
For $k=2$, Shimura's periods are defined in the proof of, and the 
remark following, \cite[Theorem 2]{shimura2}. 
It follows from the theorem above and Deligne's conjecture 
that Shimura's periods $u^{\pm}(\varphi)$ may be identified, up to 
elements of ${\mathbb Q}(\varphi)$, with Deligne's periods 
$c^{\pm}(\varphi)$.

\begin{thm}[Symmetric squared; Sturm \cite{sturm1}, \cite{sturm2}]
\label{thm:sturm}
Let $\varphi$ be a primitive form in $S_k(N,\omega)$. 
Let $\varphi(z) = \sum_{n=1}^{\infty} a_nq^n$ be its Fourier 
expansion. Let $\xi$ be a Dirichlet character. The Euler product
$$
L_f(s, {\rm Sym}^2 \varphi, \xi) = \prod_p 
\prod_{i=0}^2
(1-\alpha_{p,\varphi}^i\beta_{p,\varphi}^{2-i}\xi(p)p^{-s})^{-1}
$$
converges for ${\rm Re}(s) \gg 0$, has a meromorphic continuation and has
at most two simple poles at $s=k$ and $s=k-1$. 
Let $\nu \in \{0,1\}$ be defined by
$\xi(-1)=(-1)^{\nu}$. Let
$$
Z(s, {\rm Sym}^2 \varphi, \xi) = 
\frac{L_f(s, {\rm Sym}^2 \varphi, \xi) \pi^{k-2m-2}}
{\langle \varphi , \varphi \rangle \gamma(\theta^2)}
$$
where $\theta(a) = \omega(a)\xi(a)\left(\frac{-1}{a}\right)^{k+\nu}$
and $\gamma(\theta^2)$ is the Gauss sum associated to $\theta^2$.
For every $\sigma \in {\rm Aut}({\mathbb C})$, and 
any integer $m$ with $k \leq m \leq 2k-2-\nu$
and $m \equiv \nu \pmod{2}$, we have 
$$
Z(m, {\rm Sym}^2 \varphi, \xi)^{\sigma} = 
Z(m, {\rm Sym}^2 \varphi^{\sigma}, \xi^{\sigma}).
$$
\end{thm}

Sturm proves this result in \cite{sturm2} 
for $m = k$ and for the rest of the values of $m$ it has been 
known from his earlier paper \cite{sturm1}. In 
\cite{sturm2} there is a typo and the exponent of $\pi$ is incorrectly
written as $2m+1-k$, the correct one, which is $2m+2-k$
may be found in his earlier paper \cite{sturm1}. (Incidentally, there is 
an amusing typo of a different nature in \cite{sturm1}; see page 221 therein.)
Using (2) of Theorem~\ref{thm:shimura}, and Lemma~\ref{lem:gauss-sum} 
below, we can rewrite Sturm's theorem as 
{\Small
\begin{equation}
\label{eqn:sturm}
\sigma\left(\frac{L_f(m, {\rm Sym}^2 \varphi, \xi)}
{(2\pi i)^{2m+1-k} (u^+(\varphi)u^-(\varphi))
\gamma(\omega \xi^2)}\right) =
\left(\frac{L_f(m, {\rm Sym}^2 \varphi^{\sigma}, \xi^{\sigma})}
{(2\pi i)^{2m+1-k} (u^+(\varphi^{\sigma})u^-(\varphi^{\sigma}))
\gamma((\omega\xi^2)^{\sigma})}\right).
\end{equation}}
Observe that if $\xi$ is trivial, then using Lemma~\ref{lem:gauss-sum}, 
the above equation exactly says that Conjecture~\ref{conj:deligne} is 
true for $n=2$.

\subsection{Some lemmas}
\label{sec:lemmas}

The purpose of this subsection is to record the critical integers of a symmetric
power $L$-function associated to a modular form, as well as to record some useful
lemmas which will be needed later in the paper. 

We let $W_{\mathbb R}$ be the Weil group of ${\mathbb R}$. Recall 
that as a set it is 
defined as $W_{\mathbb R} = {\mathbb C}^* \cup j{\mathbb C}^*$. The group 
structure is induced from that of ${\mathbb C}^*$ and the relations 
$jzj^{-1} = \overline{z}$ and $j^2= -1 $. We have a homomorphism 
$W_{\mathbb R} \to {\mathbb R}^*$ which sends $z \in {\mathbb C}^*$ to 
$|z|$ and sends $j$ to $-1$. This homomorphism induces an isomorphism 
from the abelianization $W_{\mathbb R}^{\rm ab} \to {\mathbb R}^*$. 
We will let $\epsilon : W_{\mathbb R} \to \{\pm 1\}$ denote the sign 
homomorphism, defined as $\epsilon(z) = 1$ and $\epsilon(j)=-1$.
For these and other details on $W_{\mathbb R}$ we refer the reader to 
\cite{knapp}. We let $1\!\!1$ denote the trivial representation of the
group in context. 

Let $k \geq 1$ be an integer. We let $\chi_{k-1}$ denote the character of 
${\mathbb C}^*$ given by $z \mapsto (z/|z|)^{k-1}$. Let $I(\chi_{k-1})$
denote the representation 
$$
I(\chi_{k-1}) = {\rm Ind}_{{\mathbb C}^*}^{W_{\mathbb R}}(\chi_{k-1}).
$$
This is the Langlands parameter of the representation at infinity of a 
weight $k$ modular form. Observe that $I(\chi_{k-1})$ is irreducible 
if $k \geq 2$.

\begin{lemma}
\label{lem:repn-at-infty}
Let $k$ be an integer $\geq 2$. For any $n \geq 1$ we have 
\begin{enumerate}
\item If $n=2r+1$ then
$$
{\rm Sym}^n(I(\chi_{k-1})) = \bigoplus_{a=0}^r I(\chi_{(2a+1)(k-1)}).
$$ 
\item If $n=2r$ then
$$
{\rm Sym}^n(I(\chi_{k-1})) = \epsilon^{r(k-1)} \oplus
\bigoplus_{a=1}^r I(\chi_{(2a(k-1))}).
$$
\end{enumerate}
\end{lemma}

\begin{proof}
The proof is quite easy and anyway such lemmas are well known to 
experts. We sketch the details here for lack of a good reference. 
We begin by making the following observations. 
Let $\sigma = I(\chi_{(k-1)})$. The determinant of $\sigma$
is given by ${\rm det}(\sigma) = \epsilon^k$. For any two integers $a,b$, 
$I(\chi_a) \otimes I(\chi_b) \simeq I(\chi_{a+b}) \oplus I(\chi_{a-b})$
and $\epsilon \otimes I(\chi_a) \simeq I(\chi_a)$. Finally, 
$I(1\!\!1) = 1\!\!1 \oplus \epsilon$. 
For any two dimensional 
representation $\sigma$ we have
\begin{eqnarray*}
\sigma \otimes \sigma 
& \simeq & 
{\rm Sym}^2(\sigma) \oplus
{\rm det}(\sigma), \\
{\rm Sym}^2(\sigma) \otimes \sigma
& \simeq &
{\rm Sym}^3(\sigma) \oplus
(\sigma \otimes {\rm det}(\sigma)). 
\end{eqnarray*}
This proves the cases $n=2,3$. For any $n$, observe that
$$
{\rm Sym}^n(\sigma) \otimes \sigma \simeq
{\rm Sym}^{n+1}(\sigma) \oplus
{\rm Sym}^{n-1}(\sigma)\otimes {\rm det}(\sigma).$$
The proof follows by induction on $n$.
\end{proof}

\begin{lemma}[Archimedean local factors \cite{knapp}]
\label{lem:archimedean}
Let $\sigma$ be an irreducible representation of $W_{\mathbb R}$. The
local factor $L(s,\sigma)$ is given by the following formulae:
$$
L(s,\sigma) = \left\{\begin{array}{ll}
\pi^{-(s+t)/2} \, \Gamma\left(\frac{s+t}{2}\right) &
\mbox{, if $\sigma = |\cdot |_{\mathbb R}^t$} \\
\pi^{-(s+t+1)/2} \, \Gamma\left(\frac{s+t+1}{2}\right) &
\mbox{, if $\sigma = \epsilon \otimes |\cdot |_{\mathbb R}^t$} \\
2(2\pi)^{-(s+t+l/2)}\ \, \Gamma(s+t + l/2 ) &
\mbox{, if $\sigma = I(\chi_l) \otimes |\cdot |_{\mathbb R}^t$ with
$l \geq 1$.}
\end{array}\right.
$$
\end{lemma}

\begin{lemma}
\label{lem:classical-modern}
Let $\varphi \in S_k(N,\omega)$ be a primitive form and let 
$\pi(\varphi)$ be the associated cuspidal automorphic representation. 
Let $n \geq 1$ be an integer. 
We have the following equality of $L$-functions: 
$$
L(s, {\rm Sym}^n, \pi(\varphi)) = L(s + n(k-1)/2, {\rm Sym}^n \varphi). 
$$
The left hand side is the $L$-function attached by Langlands to 
$\pi(\varphi)$ corresponding to the $n$-th symmetric power 
of the standard representation of ${\rm GL}_2({\mathbb C})$.

\end{lemma}

\begin{proof}
We have 
$$
L_f(s, \varphi) = 
\prod_p L_p(s, \varphi) = 
(1-\alpha_{p,\varphi} p^{-s})^{-1}(1-\beta_{p,\varphi} p^{-s})^{-1}.
$$
Similarly, we have the $L$-function of 
$\pi = \pi(\varphi) = \otimes_{p \leq \infty}' \pi_p$ given by 
$$
L_f(s, \pi) = \prod_p L_p(s, \pi_p) = 
\prod_p (1-\alpha_{p,\pi} p^{-s})^{-1}(1-\beta_{p,\pi} p^{-s})^{-1}.
$$

We know that $L_f(s, \pi) = L_f(s + (k-1)/2, \varphi)$. 
(See \cite[Example 6.19]{gelbart1} for instance.)
Hence $\alpha_{p,\pi} = \alpha_{p,\varphi}p^{-(k-1)/2}$ and 
$\beta_{p,\pi} = \beta_{p,\varphi} p^{-(k-1)/2}$. The lemma follows from 
the Euler products for both the symmetric power $L$-functions. 
\end{proof}

We can now record the critical integers for symmetric power $L$-functions. 
The main ingredients of the proof involves 
the local Langlands correspondence at infinity and some of
the lemmas above. An artifice one keeps in mind is that for an 
$L$-function of an automorphic representation of ${\rm GL}_n$, we look for 
critical points $m$ which are integral, if $n$ is odd, and are 
half-integral of the form $m + 1/2$, if $n$ is even. This artifice corresponds
to the so-called {\it motivic normalization} \cite[p. 139]{clozel2}. 

\begin{lemma}
\label{lem:critical-sym-odd}
Let $\varphi$ be a primitive cusp form of weight $k$. The set of critical integers
for $L_f(s, {\rm Sym}^{2r+1} \varphi)$ is given by integers $m$ with
$$
r(k-1)+1 \leq m \leq (r+1)(k-1). 
$$
\end{lemma}

\begin{proof}
See the proof of Lemma~\ref{lem:critical-sym-even} below.
\end{proof}

\begin{lemma}
\label{lem:critical-sym-even} 
Let $\varphi$ be a primitive cusp form of weight $k$.
The set of critical integers for $L_f(s, {\rm Sym}^{2r} \varphi)$ is given by:
{\small
\begin{enumerate}
\item If $r$ odd and $k$ even then
$$
\{(r-1)(k-1)+1, (r-1)(k-1)+3,\dots, r(k-1);\  r(k-1)+1, r(k-1)+3,\dots, (r+1)(k-1)\}.
$$
\item If $r$ odd and $k$ odd then
$$
\{(r-1)(k-1)+1, (r-1)(k-1)+3,\dots, r(k-1)-1;\  r(k-1)+2, r(k-1)+4,\dots, (r+1)(k-1)\}.
$$
\item If $r$ even and $k$ even then
$$
\{(r-1)(k-1)+2, (r-1)(k-1)+4,\dots, r(k-1)-1;\  r(k-1)+2, r(k-1)+4,\dots, (r+1)(k-1)-1\}.
$$
\item If $r$ even and $k$ odd then
$$
\{(r-1)(k-1)+1, (r-1)(k-1)+3,\dots, r(k-1)-1;\  r(k-1)+2, r(k-1)+4,\dots, (r+1)(k-1)\}.
$$
\end{enumerate}}
\end{lemma}

\begin{proof}

The proof is a rather tedious application of the above lemmas. As a 
representative example we prove it for $L(s, {\rm Sym}^4 \varphi)$ and 
leave the general case to the reader!
Let $\pi = \pi(\varphi)$ be the representation associated to $\varphi$. 
We will identify the critical points for $L(s, {\rm Sym}^4, \pi)$,
and the corresponding statement 
for $L(s, {\rm Sym}^4 \phi)$ follows from 
Lemma~\ref{lem:classical-modern}. By the above mentioned `artifice' we 
look for integers $m$ such that 
$L_{\infty}(s, {\rm Sym}^4, \pi)$ and 
$L_{\infty}(1-s, {\rm Sym}^4, \pi^{\vee})$ are regular 
at $s=m$. 

Via the local Langlands correspondence 
we transfer our attention to the $L$-functions at 
infinity on the `Galois side'. Since $\pi = \pi(\varphi)$, the 
representation $\pi_{\infty}$ is a discrete series representation of 
${\rm GL}_2({\mathbb R})$ of lowest weight $k$. The Langlands parameter of this is the 
representation 
$I(\chi_{(k-1)}) = 
{\rm Ind}_{{\mathbb C}^*}^{W_{\mathbb R}}(\chi_{k-1})$. 
Hence 
$$
L_{\infty}(s, {\rm Sym}^4, \pi) = 
L(s, {\rm Sym}^4(I(\chi_{(k-1)}))).
$$
Using Lemma~\ref{lem:repn-at-infty} and Lemma~\ref{lem:archimedean}, we 
get 
$$
L(s, {\rm Sym}^4(I(\chi_{(k-1)}))) 
\sim
\Gamma(s+2(k-1))\Gamma(s+k-1)\Gamma(s/2).
$$
Just for this proof, by $\sim$ we mean up to an exponential function, 
which is holomorphic and nonvanishing everywhere, and 
so is irrelevant for the computation of a critical point. 
We also have 
$$
L_{\infty}(1-s, {\rm Sym}^4, \pi^{\vee}) 
\sim \Gamma(2k-s-1)\Gamma(k-s)\Gamma((1-s)/2).
$$
since the Langlands correspondence and symmetric powers both commute with 
taking contragredients. Hence, we get that an integer $m$ is critical if 
\begin{enumerate}
\item $m+k-1 \geq 1$,
\item $m$ is not an even nonpositive integer, 
\item $k-m \geq 1$, and 
\item $m$ not an odd positive integer. 
\end{enumerate}
\end{proof}

We end this subsection by 
recalling some standard facts about Dirichlet characters. 
We will identify Dirichlet characters with characters of 
the id\`ele class group of ${\mathbb Q}$ via 
the isomorphism \cite[Proposition 6.1.10]{neukirch}. 
An important detail in this dictionary is that the parity
of a Dirichlet character $\chi$ 
is seen by the infinity component $\chi_{\infty}$ of the 
corresponding id\`ele class character, i.e., $\chi(-1) = \chi_{\infty}(-1)$.
Let $\varphi$ be a primitive form in $S_k(N,\omega)$, where $\omega$ is a 
Dirichlet character modulo $N$. Let 
$\pi(\varphi)$ be the associated cuspidal automorphic representation of 
${\rm GL}_2({\mathbb A}_{\mathbb Q})$. Let $\omega_{\pi(\varphi)}$ be 
the central character of $\pi(\varphi)$; it is an id\`ele class character. 
Under the above identification, we have 
$\omega_{\pi(\varphi)} = \omega$. 
This may be seen by comparing the coeffcients of $p^{-2s}$ in the Euler products appearing
in the proof of Lemma~\ref{lem:classical-modern}. We finally recall an important 
property of the Gauss sum of the product of two Dirichlet characters.

\begin{lemma}
\label{lem:gauss-sum}
Let $\omega_1$ and $\omega_2$ be Dirichlet characters. 
For all $\sigma \in {\rm Aut}({\mathbb C})$, we have
$$
\sigma\left(\frac{\gamma(\omega_1)\gamma(\omega_2)}
{\gamma(\omega_1\omega_2)}\right) = 
\frac{\gamma(\omega_1^{\sigma})\gamma(\omega_2^{\sigma})}
{\gamma(\omega_1^{\sigma}\omega_2^{\sigma})}.
$$ 
\end{lemma}

\begin{proof}
This lemma is due to Shimura. See \cite[Lemma 8]{shimura1} where the 
proof is unreasonably complicated. A simpler proof is suggested in 
\cite[(\S4, Remark 1)]{shimura2}. For the sake of completeness we 
sketch this proof. 
Let $c$ be the least common multiple of the conductors of $\omega_1$, 
$\omega_2$ and $\omega_1\omega_2$. Let $b \in {\mathbb Z}$ be relatively 
prime to $c$ such that $\sigma(e^{2\pi i/c}) = e^{2\pi i b/c}$.
If $c_1$ is the conductor of $\omega_1$, then we have
\begin{eqnarray*}
\sigma(\gamma(\omega_1)) 
& = & \sigma(\sum_{x=0}^{c_1-1}\omega_1(x) e^{2\pi ix/c_1}) 
=  \sum_{x=0}^{c_1-1}\sigma(\omega_1(x)) e^{2\pi ibx/c_1} \\
& = & \sigma(\omega_1(b)^{-1})\sum_{x=0}^{c_1-1}\omega_1^{\sigma}(x) 
e^{2\pi ix/c_1} = 
\sigma(\omega_1(b)^{-1})\gamma(\omega_1^{\sigma}).
\end{eqnarray*}
Similarly, $\sigma(\gamma(\omega_2)) = 
\sigma(\omega_2(b)^{-1})\gamma(\omega_2^{\sigma})$ and 
$\sigma(\gamma(\omega_1\omega_2)) = 
\sigma(\omega_1(b)^{-1}\omega_2(b)^{-1})
\gamma(\omega_1^{\sigma}\omega_2^{\sigma}).$ 
Hence the quotient
$\gamma(\omega_1)\gamma(\omega_2)/\gamma(\omega_1\omega_2)$ is 
equivariant under all $\sigma \in {\rm Aut}({\mathbb C}).$
\end{proof}

\subsection{Dihedral forms}
\label{sec:dihedral}

In this section we prove a theorem about the special values of any 
symmetric power $L$-function associated to a dihedral cusp form. This 
formally looks exactly like Deligne's conjecture, the only {\it 
difference} being that Deligne's motivcally defined periods 
$c^{\pm}(\varphi)$ are replaced by Shimura's periods $u^{\pm}(\varphi)$. 
(See the paragraph after Theorem~\ref{thm:shimura}.)
As mentioned earlier, Deligne's conjecture for dihedral 
forms is known via motivic considerations. In what follows we use only
$L$-functions, and in the process use some nonvanishing 
results for twists of $L$-functions. The technical heart of the proof 
below is a certain {\it period relation} which is interesting in its own
right, and it is this relation which justifies this section. 
If $\varphi_{\chi}$ denotes the dihedral modular form corresponding to a 
character $\chi$ of an imaginary quadratic number field, then the main
theorem proved in this section relates the periods of 
$\varphi_{\chi^n}$--for any power 
$\chi^n$ of $\chi$--to the periods of $\varphi_{\chi}$
(see Theorem~\ref{thm:period-relations} below).

Note that if one has a cuspidal automorphic representation $\pi$ of 
${\rm GL}_2({\mathbb A}_{\mathbb Q})$, and suppose 
$\pi = {\rm AI}_{K/{\mathbb Q}}(\chi)$ is dihedral, then every 
symmetric power lifting ${\rm Sym}^n(\pi)$, in the sense of Langlands
functoriality, exists.
To state functoriality, we need some notation. Let $W_{{\mathbb Q}_p}'$ 
denote the Weil-Deligne group of ${\mathbb Q}_p$. The local Langlands 
correspondence (see \cite{kudla}) gives a bijection between irreducible
admissible representations $\pi$ of ${\rm GL}_n({\mathbb Q}_p)$ and 
$n$-dimensional semisimple representations $\sigma$ of 
$W_{{\mathbb Q}_p}'$. We will denote this bijection by $\pi \mapsto 
\sigma(\pi)$ and similarly $\sigma \mapsto \pi(\sigma)$. 
Now let $\pi$ denote a cuspidal automorphic representation of 
${\rm GL}_2({\mathbb A}_{\mathbb Q})$. Then $\pi$ is a tensor product 
of local representations as $\pi = \otimes' \pi_p$. (We let $p$ run 
through the finite primes as well as $\infty$.) To each $\pi_p$, a 
representation of ${\rm GL}_2({\mathbb Q}_p)$, is associated via the local 
Langlands correspondence, a representation 
$\pi({\rm Sym}^n(\sigma(\pi_p)))$ of ${\rm GL}_{n+1}({\mathbb Q}_p)$.
We will denote this representation as ${\rm Sym}^n(\pi_p)$. 
If $\pi_p$ is unramified then so is ${\rm Sym}^n(\pi_p)$.
The global symmetric power lift is defined as 
$
{\rm Sym}^n(\pi) = \otimes'_{p \leq \infty} {\rm Sym}^n(\pi_p). 
$
Langlands functoriality takes the form that the irreducible representation
${\rm Sym}^n(\pi)$ is an isobaric automorphic representation of 
${\rm GL}_{n+1}({\mathbb A}_{\mathbb Q})$.
This functorial formalism also asserts
that the symmetric power $L$-function $L(s, {\rm Sym}^n, \pi)$ is the 
standard $L$-function $L(s, {\rm Sym}^n(\pi))$ of  
${\rm Sym}^n(\pi)$.  If $\pi$ is dihedral then one can indeed
write down the isobaric decomposition of ${\rm Sym}^n(\pi)$. 
See Lemma~\ref{lem:isobaric-dihedral} below.

For the rest of this section we let $K/{\mathbb Q}$ be an imaginary 
quadratic extension. 
We will let $\omega_K = \omega_{K/{\mathbb Q}}$ denote the 
corresponding quadratic character of 
${\mathbb Q}^*\backslash {\mathbb I}_{\mathbb Q}$. 
(Note that as a Dirichlet character, $\omega_K$ is an odd character.)
We let $\gamma_K$ denote the Gauss sum associated to $\omega_K$.
We let $\chi$ denote a character of 
$K^*\backslash {\mathbb I}_K$ such that its infinity component 
is $\chi_{\infty}(z) = (z/|z|)^{k-1}$ for an integer $k \geq 2$. 
Hence $\chi^n$ is not ${\rm Gal}(K/{\mathbb Q})$ invariant for any $n \geq 
1$ by the following lemma. 

\begin{lemma}
\label{lem:weight1}
Let $\varphi$ be a primitive form in $S_k(N,\omega)$. Suppose that 
$\pi = \pi(\varphi)$ is a dihedral form, 
$\pi = {\rm AI}_{K/{\mathbb Q}}(\chi)$, where for some integer $r \geq 1$,
$\chi^r$ is invariant under the Galois group of $K/{\mathbb Q}$. 
Then $\varphi$ is necessarily a weight $1$ form, i.e., $k=1$. 
\end{lemma}

\begin{proof}
Since $\pi = \pi(\varphi) = {\rm AI}_{K/{\mathbb Q}}(\chi)$
corresponds to a weight $k$ form, we must have the
character at infinity $\chi_{\infty} : K_{\infty} \to {\mathbb C}^*$, 
to be given by $\chi_{\infty}(z) = (z/|z|)^{k-1}.$
Since $\chi^r$ is Galois invariant, there is an id\`ele class character $\mu$ of 
${\mathbb Q}$ such that $\chi^r = \mu \circ N_{K/{\mathbb Q}}$. 
Let the character $\mu_{\infty} : {\mathbb Q}_{\infty} \to {\mathbb C}^*$ 
be given by $\mu_{\infty}(a) = {\rm sign}(a)^{\varepsilon}|a|^t$
for $\varepsilon \in \{0,1\}$ and some $t \in {\mathbb C}$. 
The relation 
$\chi_{\infty}^r = \mu_{\infty}\circ 
N_{K_{\infty}/{\mathbb Q}_{\infty}}$ gives 
$(z/|z|)^{r(k-1)} = |z|^{2t}$ for all $z \in {\mathbb C}^*$.
It is easy to see that this forces $k=1$ and $t=0$. 
\end{proof}

We note that as far as Deligne's conjectures are concerned, weight $1$ forms are 
not interesting since by 
Lemma~\ref{lem:critical-sym-odd} and Lemma~\ref{lem:critical-sym-even}
the symmetric power $L$-functions do not have 
critical points. 
We will henceforth assume that 
$\chi^n$ is not Galois invariant for any nonzero integer $n$.

Let $\pi = \pi(\chi)  := {\rm AI}_{K/{\mathbb Q}}(\chi)$ 
be the dihedral cuspidal automorphic representation associated to $\chi$. 
We denote the corresponding 
holomorphic cusp form as $\varphi_{\chi}$. 
Note that $\varphi_{\chi} \in S_k(N,\omega)$ where the level $N$ depends 
on the conductor of $\chi$ and the discriminant of $K$, and 
the nebentypus $\omega$ 
can be described as $\omega \omega_K = 
\chi_{\mathbb Q}$ as an equality of id\`ele class characters of 
${\mathbb Q}$. (Here $\chi_{\mathbb Q}$ denotes the restriction of 
$\chi$ to ${\mathbb I}_{\mathbb Q}$.) 
The first step to proving Deligne's conjecture 
is to write down a decomposition of $L(s, {\rm Sym}^n \varphi)$
when $\varphi = \varphi_{\chi}$ is dihedral.
The following lemma gives the isobaric decomposition of a symmetric power lift of a 
dihedral cusp form. 

\begin{lemma}
\label{lem:isobaric-dihedral}
Let $\chi$ be a character of a quadratic extension $K/{\mathbb Q}$ and assume that
$\chi^n$ is not Galois invariant for any nonzero integer $n$. Recall that 
$\chi_{\mathbb Q}$
is the restriction of $\chi$ to the id\`eles of ${\mathbb Q}$. 
Then we have
\begin{eqnarray*}
{\rm Sym}^{2r}({\rm AI}_{K/{\mathbb Q}}(\chi)) & = & 
\boxplus_{a=0}^{r-1}{\rm AI}_{K/{\mathbb Q}}(\chi^{2r-a}\chi'^a) 
\boxplus \chi_{\mathbb Q}^r, \\
{\rm Sym}^{2r+1}({\rm AI}_{K/{\mathbb Q}}(\chi)) & = &
\boxplus_{a=0}^{r}{\rm AI}_{K/{\mathbb Q}}(\chi^{2r+1-a}\chi'^a).
\end{eqnarray*}
\end{lemma}

\begin{proof}
The proof is by induction on $n$ for ${\rm Sym}^n({\rm AI}_{K/{\mathbb 
Q}}(\chi))$ and 
is analogous to the proof of Lemma~\ref{lem:repn-at-infty}. 
We leave the details to the reader.
\end{proof}

\begin{cor}
\label{cor:dihedral}
The symmetric power $L$-functions of a dihedral cusp form decompose as 
follows: 
\small{
\begin{eqnarray*}
L_f(s, {\rm Sym}^{2r} \varphi_{\chi}) 
& = & 
L_f(s-r(k-1), (\omega\omega_K)^r)
\prod_{a=0}^{r-1} L_f(s - a(k-1), \varphi_{\chi^{2(r-a)}}, \omega^a) \\
& = & 
L_f(s-r(k-1), (\omega\omega_K)^r)
\prod_{a=0}^{r-1} L_f(s - a(k-1), \varphi_{\chi^{2(r-a)}}, 
(\omega \omega_K)^a). \\
L_f(s, {\rm Sym}^{2r+1} \varphi_{\chi})
& = &
\prod_{a=0}^{r} L_f(s - a(k-1), \varphi_{\chi^{2(r-a)+1}}, \omega^a) \\
& = &
\prod_{a=0}^{r} L_f(s - a(k-1), \varphi_{\chi^{2(r-a)+1}}, 
(\omega \omega_K)^a).
\end{eqnarray*}
}
\end{cor}

\begin{proof}
Note that
${\rm AI}_{K/{\mathbb Q}}(\chi^{2r-a}\chi'^a) \simeq
{\rm AI}_{K/{\mathbb Q}}(\chi^{2(r-a)}) \otimes \chi_{\mathbb Q}^a$ 
and a similar statement for odd symmetric powers. Note also that
for any integer $l$,
${\rm AI}_{K/{\mathbb Q}}(\chi^l)
\simeq {\rm AI}_{K/{\mathbb Q}}(\chi^l) \otimes \omega_K$.
The proof follows from Lemma~\ref{lem:classical-modern} and 
Lemma~\ref{lem:isobaric-dihedral}. 
\end{proof}

To prove Deligne's conjecture for dihedral forms, 
we will need to relate the periods of the 
cusp forms $\varphi_{\chi^n}$ to the periods of the cusp form 
$\varphi_{\chi}$, but before doing so we need some preliminaries on 
Galois properties of dihedral forms. Especially, we want to know 
the behaviour of dihedral forms under the action of
${\rm Aut}({\mathbb C})$.

Given a primitive modular form $\varphi$ of weight $k$, with Fourier expansion
$\varphi(z) = \sum a_n e^{2\pi i nz}$, and given 
$\sigma \in {\rm Aut}({\mathbb C})$, we define
$\varphi^{\sigma}(z) = \sum \sigma(a_n)e^{2\pi i nz}$, which is also a
primitive modular form of the same weight $k$. 
We begin by observing that the process of attaching
a cuspidal representation $\pi(\varphi)$ to $\varphi$ is not an equivariant
process in general. It is so exactly when the weight $k$ is even, which is 
also the parity condition on $k$ which ensures that $\pi(\varphi)$ is 
algebraic. It is easily checked that $\pi(\varphi)\otimes |\!|\ |\!|^{-k/2}$ is 
an algebraic (regular cuspidal automorphic) representation. (See Clozel \cite[p.91]{clozel2}.) Appealing to \cite[Th\'eor\`eme 3.13]{clozel2} we deduce that 
$(\pi(\varphi)\otimes |\!|\ |\!|^{-k/2})^{\sigma}$ 
is an algebraic cuspidal representation. Indeed, we have
\begin{equation}
\label{eqn:sigma-repn}
(\pi(\varphi)\otimes |\!|\ |\!|^{-k/2})^{\sigma}
\ = \ 
\pi(\varphi^{\sigma})\otimes |\!|\ |\!|^{-k/2}.
\end{equation}
This may be seen by comparing both sides at all unramified places, 
while using Waldspurger \cite[Exemple \S I.2]{waldspurger1}. 
(It is interesting to note that, in the spirit of Clozel \cite[Definitions 
1.9--1.11]{clozel2}, one can define the process $\varphi \mapsto \pi(\varphi)$, 
with a Tate twist  
$\pi^T(\varphi) := \pi(\varphi)\otimes |\!|\ |\!|^{-k/2}$, so that the map
$\varphi \mapsto \pi^T(\varphi)$ is equivariant, i.e., respects algebraicity.)

Next, we analyze such an equivariance property for automorphic induction.  
Let $\chi$ be a Hecke character of $K$ (an imaginary quadratic extension) with 
$\chi_{\infty}(z) = (z/|z|)^{k-1}$ for an integer $k \geq 2$. Consider the
automorphic induction ${\rm AI}_{K/{\mathbb Q}}(\chi)$, which is a cuspidal
representation. As above, ${\rm AI}_{K/{\mathbb Q}}(\chi)\otimes |\!|\ |\!|^{-k/2}$
is an algebraic representation. We can apply $\sigma$ to this, and ask for 
the relation of the resulting representation with the induction of $\chi^\sigma$. 
Note that $\chi$ is not algebraic in general, however, 
$\chi \otimes |\!|\ |\!|^{-(k-1)/2}$ is an algebraic id\`ele class character, and so we can
apply $\sigma$ to such a twist of $\chi$. 
(For the definition of $(\chi \otimes |\!|\ |\!|^{-(k-1)/2})^{\sigma}$ see
Clozel \cite[p.107]{clozel2}.) We have
\begin{equation}
\label{eqn:sigma-ai}
({\rm AI}_{K/{\mathbb Q}}(\chi)\otimes |\!|\ |\!|^{-k/2})^{\sigma} = 
{\rm AI}_{K/{\mathbb Q}}((\chi \otimes |\!|\ |\!|^{-(k-1)/2})^{\sigma})\otimes |\!|\ |\!|^{-1/2}.
\end{equation}

To such a character $\chi$ we have the modular cusp form $\varphi_{\chi}$, which 
we recall is defined as that form for which 
$\pi(\varphi_{\chi}) = {\rm AI}_{K/{\mathbb Q}}(\chi)$. 
From (\ref{eqn:sigma-repn}) and (\ref{eqn:sigma-ai}) we deduce 
\begin{equation}
\label{eqn:sigma-form}
\varphi_{\chi}^{\sigma} = 
\varphi_{(\chi\otimes |\!|\ |\!|^{-(k-1)/2})^{\sigma}\otimes |\!|\ |\!|^{(k-1)/2}}.
\end{equation}
{\it In particular, if $\varphi$ is a dihedral form then so is $\varphi^{\sigma}$, and comes
from the same quadratic extension (namely $K$) as that for $\varphi$.} 

(If $k$ is odd, then (\ref{eqn:sigma-form}) simplifies to 
$\varphi_{\chi}^{\sigma} = \varphi_{\chi^{\sigma}}$. 
This can be seen easily in the classical setup: We let 
$\tilde{\chi}$ be the corresponding Gr\"ossencharakter attached to $\chi$. 
Define a function $\varphi_{\tilde{\chi}}$ on the
upper half plane by 
$$
\varphi_{\tilde{\chi}}(z) = 
\sum_{\mathfrak{a}}
\tilde{\chi}(\mathfrak{a})
N(\mathfrak{a})^{(k-1)/2}e^{2\pi i N(\mathfrak{a})z}
$$
where $\mathfrak{a}$ runs over all the integral ideals of $K$.
Assume that $\chi$ (or equivalently $\tilde{\chi}$) is primitive. Then
$\varphi_{\tilde{\chi}}$ is a primitive modular cusp form. 
(See \cite[Theorem 3.8.2]{miyake}.) Then 
$\varphi_{\chi} = \varphi_{\tilde{\chi}}$ which may be seen by 
comparing Satake parameters for both the modular forms. Now the equivariance
of $\chi \mapsto \varphi_{\chi}$ is obvious when $k$ is odd.)

\smallskip

We can now state and prove the main result of this section. 

\begin{thm}[Period relations for dihedral forms]
\label{thm:period-relations}
Let $\chi$ be a Hecke character of an imaginary quadratic field $K$ with 
$\chi_{\infty}(z) = (z/|z|)^{k-1}$ for an integer $k \geq 2$. Let $\varphi_{\chi}$
be the corresponding modular cusp form. For any positive integer $n$ and for all
$\sigma \in {\rm Aut}({\mathbb C})$ we have 
\begin{enumerate}
\item $$
\sigma\left(\frac{u^+(\varphi_{\chi^n})}{u^+(\varphi_{\chi})^n}\right) = 
\frac{u^+(\varphi^{\sigma}_{\chi^n})}{u^+(\varphi^{\sigma}_{\chi})^n} $$
\item $$
\sigma\left(\frac{u^-(\varphi_{\chi^n})}{u^+(\varphi_{\chi})^n \gamma_K}\right) = 
\frac{u^-(\varphi^{\sigma}_{\chi^n})}{u^+(\varphi^{\sigma}_{\chi})^n\gamma_K} $$
\end{enumerate}
where $\gamma_K$ is the Gauss sum of the quadratic character $\omega_K$ of 
${\mathbb Q}$ associated to $K$.
\end{thm}

\begin{proof}
Before we get into the proof of the theorem, we record some nonvanishing results 
for $L$-functions which will be useful later on. 
\begin{lemma}
\label{lem:nonvanishing}
With the notations as above
\begin{enumerate}
\item There is an even Dirichlet character $\xi$ 
such that $L_f(1,\varphi_{\chi}, \xi) \neq 0$.
\item There is an even Dirichlet character $\xi$ such that  
$L_f(k, \varphi_{\chi^n} \times (\varphi_{\chi} \otimes \xi)) \neq 0.$
\end{enumerate}
\end{lemma}

\begin{proof}
The first assertion follows by using the main theorem of Rohrlich 
\cite{rohrlich}. The second assertion follows by using the main theorem 
of Barthel--Ramakrishnan \cite{barthel-ramakrishnan} while thinking of the 
Rankin--Selberg $L$-function as a standard $L$-function for ${\rm GL}_4$, which 
we can do by the work of Ramakrishnan \cite{ramakrishnan}.
\end{proof}

We will prove the theorem by induction on $n$. The following lemma is the 
$n=1$ case of the theorem, for which (1) is a tautology, and it is only 
statement (2) which needs a proof. Observe that (2) implies in particular 
that for any dihedral form the periods $u^{\pm}$ are algebraically 
dependent. (See also Harris \cite[Remark (2.7)]{harris} and 
Bertrand \cite[Corollary 1, p.35]{bertrand}.)

\begin{lemma}
\label{lem:dependence}
With the notations as in Theorem~\ref{thm:period-relations} we have
$$
\sigma\left(\frac{u^{\pm}(\varphi_{\chi})}{u^{\mp}(\varphi_{\chi}) \gamma_K}\right) =
\frac{u^{\pm}(\varphi^{\sigma}_{\chi})}{u^{\mp}(\varphi^{\sigma}_{\chi})\gamma_K}.
$$
\end{lemma}

\begin{proof}
By Lemma~\ref{lem:nonvanishing} there is 
an even Dirichlet character $\xi$ such that $L(1, \varphi_{\chi}, \xi) \neq 0.$
Since ${\rm AI}_{K/{\mathbb Q}}(\chi) \simeq {\rm AI}_{K/{\mathbb Q}}(\chi) \otimes 
\omega_K$ we get $L_f(s, \varphi_{\chi}, \eta) = L_f(s, \varphi_{\chi}, \eta\omega_K)$
for any Dirichlet character $\eta$. We have 
$$
\frac{u^+(\varphi_{\chi})\gamma_K}{u^-(\varphi_{\chi})} = 
\left(\frac{(2\pi i)u^+(\varphi_{\chi}) \gamma(\xi \omega_K)}
{L(1, \varphi_{\chi}, \xi\omega_K)}\right) 
\left(\frac{L(1, \varphi_{\chi}, \xi)}
{(2\pi i)u^-(\varphi_{\chi}) \gamma(\xi)}\right)
\left(\frac{\gamma(\xi)\gamma(\omega_K)}{\gamma(\xi\omega_K)}\right).
$$
By Theorem~\ref{thm:shimura} and Lemma~\ref{lem:gauss-sum} every factor 
on the right is ${\rm Aut}({\mathbb C})$-equivariant, and hence so is the left 
hand side. Observe also that
$$
\frac{u^+(\varphi_{\chi})}{u^-(\varphi_{\chi})\gamma_K} = 
\left(\frac{u^+(\varphi_{\chi})\gamma_K}{u^-(\varphi_{\chi})}\right)
\gamma_K^{-2}.
$$
This proves the lemma since $\gamma_K^2 \in {\mathbb Q}$.
\end{proof}

The above lemma applied to $\chi^n$ 
gives $(1) \Rightarrow (2)$ of Theorem~\ref{thm:period-relations} since
$$
\frac{u^-(\varphi_{\chi^n})}{u^+(\varphi_{\chi})^n \gamma_K} = 
\left(\frac{u^-(\varphi_{\chi^n})}{u^+(\varphi_{\chi^n}) \gamma_K}\right)
\left(\frac{u^+(\varphi_{\chi^n})}{u^+(\varphi_{\chi})^n}\right).
$$

It remains to prove (1), which we do so by 
induction on $n$. We just verified the $n=1$ case in the above lemma. Next, we  
prove it for $n=2$, and then prove it by induction for all $n$ (since the statement for
$n+1$ will depend on the statements for $n$ and $n-1$).
Applying Corollary~\ref{cor:dihedral} we have
$$
L_f(2k-2, \varphi_{\chi^2})L_f(k-1, \omega \omega_K) =
L_f(2k-2, {\rm Sym}^2 \varphi_{\chi}).
$$
Observe that $(\omega \omega_K)(-1) = (-1)^{k+1}$ and hence
$2k-2$ is critical for $L_f(s, {\rm Sym}^2 \varphi_{\chi})$ by Theorem~\ref{thm:sturm}
or by Lemma~\ref{lem:critical-sym-even}, and (necessarily) 
$k-1$ is critical for $L_f(s,\omega \omega_{K/{\mathbb Q}})$ and $2k-2$ is critical
for $L_f(s, \varphi_{\chi^2})$. Observe also that since $k \geq 2$ both the factors
on the left hand side are nonzero. Finally, in the expression
\begin{eqnarray*}
\frac{u^+(\varphi_{\chi^2})}{u^+(\varphi_{\chi})^2} 
& = &
  \left(\frac{(2\pi i)^{2k-2}u^+(\varphi_{\chi^2})}{L_f(2k-2,\varphi_{\chi^2})}\right) 
\cdot
  \left(\frac{L_f(2k-2, {\rm Sym}^2 \varphi_{\chi})}
             {(2\pi i)^{3k-3}u^+(\varphi_{\chi})u^-(\varphi_{\chi})\gamma(\omega)}\right)
\cdot \\ & & 
  \left(\frac{(2\pi i)^{k-1}\gamma(\omega\omega_K)}{L_f(k-1,\omega\omega_K)}\right) \cdot
  \left(\frac{\gamma(\omega)\gamma(\omega_K)}{\gamma(\omega\omega_K)}\right)\cdot
  \left(\frac{u^-(\varphi_{\chi})}{u^+(\varphi_{\chi})\gamma_K}\right)
\end{eqnarray*}
applying Theorem~\ref{thm:shimura}, Theorem~\ref{thm:sturm}, (\ref{eqn:dirichlet}),
Lemma~\ref{lem:gauss-sum}, and Lemma~\ref{lem:dependence}, we see that 
each of the five factors on the right is ${\rm Aut}({\mathbb C})$-equivariant, and hence 
so is the left hand side.
To apply induction for $n \geq 3$, we need the following lemma. 

\begin{lemma}
\label{lem:period-relation-induction}
For all $n \geq 2$ we have 
$$
\sigma\left(
\frac{u^+(\varphi_{\chi^{n+1}})u^-(\varphi_{\chi^{n-1}})}
{u^+(\varphi_{\chi^n})u^-(\varphi_{\chi^n})}\right) = 
\frac{u^+(\varphi^{\sigma}_{\chi^{n+1}})u^-(\varphi^{\sigma}_{\chi^{n-1}})}
{u^+(\varphi^{\sigma}_{\chi^n})u^-(\varphi^{\sigma}_{\chi^n})}
$$
\end{lemma}
 
\begin{proof} Consider the Rankin-Selberg $L$-function $L(s, 
\varphi_{\chi^n} \times \varphi_{\chi})$. It is easy to see using 
Lemma~\ref{lem:classical-modern} that 
\begin{equation} 
\label{eqn:rankin-selberg} 
L_f(s, \varphi_{\chi^n} \times \varphi_{\chi}) = 
L_f(s, \varphi_{\chi^{n+1}}) L_f(s-k+1, \varphi_{\chi^{n-1}}, \omega). 
\end{equation} 
Note that $s=k$ is critical for $L_f(s, \varphi_{\chi^{n+1}})$ and 
$L_f(s-k+1, \varphi_{\chi^{n-1}}, \omega)$ and hence it is critical
for the Rankin-Selberg $L$-function also. 
The lemma {\it follows} 
by evaluating (\ref{eqn:rankin-selberg}) at $s=k$. 
Applying Lemma~\ref{lem:nonvanishing} we can choose  
an even Dirichlet character $\xi$ such that 
$L_f(k, \varphi_{\chi^n} \times (\varphi_{\chi} \otimes \xi)) \neq 0$.
If $\pm = (-1)^k$, we have 
\begin{eqnarray*}
\frac{u^+(\varphi_{\chi^n})u^-(\varphi_{\chi^n})}
{u^{\pm}(\varphi_{\chi^{n+1}})u^{\mp}(\varphi_{\chi^{n-1}})} 
& = & 
\left(\frac{(2\pi i)^{k+1}(u^+(\varphi_{\chi^n})u^-(\varphi_{\chi^n}))
             \gamma(\omega \xi^2)}
           {L_f(k, \varphi_{\chi^n} \times (\varphi_{\chi} \otimes \xi))}\right) \cdot
\\ 
& & \left(\frac{L_f(k,\varphi_{\chi^{n+1}}, \xi)}
{(2\pi i)^k u^{\pm}(\varphi_{\chi^{n+1}}) \gamma(\xi)}\right)\cdot 
\left(\frac{L_f(1,\varphi_{\chi^{n-1}}, \omega\xi)}
{(2\pi i) u^{\mp}(\varphi_{\chi^{n-1}}) \gamma(\omega\xi)}\right)\cdot \\
& & 
\left(\frac{\gamma(\omega\xi)\gamma(\xi)}{\gamma(\omega \xi^2)}\right).
\end{eqnarray*}
In the right hand side, we see that the first factor is 
${\rm Aut}({\mathbb C})$-equivariant by applying
Shimura~\cite[Theorem 4]{shimura2}. (In the notations of that theorem, take 
$f=\varphi_{\chi^n}$, and $g = \varphi_{\chi}\otimes \xi$; observe that all the 
hypothesis of that theorem are indeed satisfied.) Further, applying 
Theorem~\ref{thm:shimura} and Lemma~\ref{lem:gauss-sum} we see that every factor
is equivariant and hence so is the left hand side. If $k$ is even then this exactly
proves the lemma, and if $k$ is odd, we still get the lemma by appealing (twice) to 
Lemma~\ref{lem:dependence}. 
\end{proof}

Theorem~\ref{thm:period-relations} follows by induction using 
Lemma~\ref{lem:period-relation-induction} (and Lemma~\ref{lem:dependence}). 
\end{proof}

We are now in a position to verify Deligne's conjecture for a dihedral cusp 
form. 

\begin{thm}
\label{thm:dihedral-deligne}
For a dihedral cusp form (like $\varphi_{\chi}$) 
Conjecture~\ref{conj:deligne} is true. 
\end{thm}

\begin{proof}
The proof follows from Corollary~\ref{cor:dihedral}, Lemma~\ref{lem:critical-sym-odd}, 
Lemma~\ref{lem:critical-sym-even},  
Theorem~\ref{thm:period-relations} and Theorem~\ref{thm:shimura}. While using 
Corollary~\ref{cor:dihedral} it is convenient to use the decompositions 
with $\omega^a$ if $k$ is even and to use those with $(\omega\omega_K)^a$ if 
$k$ is odd, since the special values of a twisted $L$-function are easy to describe 
if the twisting character is even. Carrying out the proof is rather tedious; we just 
sketch the details in one case, the rest of the cases being absolutely similar. 

Consider $L_f(m, {\rm Sym}^{2r+1} \varphi_{\chi})$. Checking the details will depend
on the parities of $m,r$ and $k$; eight cases in all. We sketch the details when 
all of them are even. 
From Corollary~\ref{cor:dihedral} we have
$$
L_f(m, {\rm Sym}^{2r+1} \varphi_{\chi}) = 
\prod_{a=0}^{r} L_f(m - a(k-1), \varphi_{\chi^{2(r-a)+1}}, \omega^a)
$$
where $m$ is given by Lemma~\ref{lem:critical-sym-odd}. 
Applying Theorem~\ref{thm:shimura} to every factor on the right we get that
$$
\prod_{a=0}^r \frac{L_f(m - a(k-1), \varphi_{\chi^{2(r-a)+1}}, \omega^a)}
{(2\pi i)^{m-a(k-1)}u^{(-1)^a}(\varphi_{\chi^{2(r-a)+1}}) \gamma(\omega^a)}
$$
is ${\rm Aut}({\mathbb C})$-equivariant. 
The denominator, after grouping together the various powers 
of $(2 \pi i)$, the powers of $\gamma(\omega)$ (using Lemma~\ref{lem:gauss-sum})
and finally the periods $u^{\pm}$, is up to equivariant quantities, the same as  
$$
(2\pi i)^{m(r+1)-(k-1)r(r+1)/2} \gamma(\omega)^{r(r+1)/2} 
(u^+(\varphi_{\chi^{2r+1}})u^-(\varphi_{\chi^{2r-1}})\cdots 
u^-(\varphi_{\chi^3})u^+(\varphi_{\chi})). 
$$
Using Theorem~\ref{thm:period-relations} and the definition of $\delta(\omega)$ this 
simplifies, up to equivariant quantities, to 
$$
(2\pi i)^{m(r+1)} \delta(\omega)^{r(r+1)/2} u^+(\varphi_{\chi})^{(r+1)^2} 
\gamma_K^{r/2}.
$$
Since $r$ is even, $r/2 \equiv r(r+1)/2 \pmod{2}$, and $\gamma_K^2 \in {\mathbb 
Q}$, using Theorem~\ref{thm:period-relations} this further simplifies to 
$$
(2\pi i)^{m(r+1)} c^+({\rm Sym}^{2r+1} \varphi_{\chi})
$$
which concludes the proof in this case. (In the definition of 
$c^+({\rm Sym}^{2r+1} \varphi_{\chi})$, we replace $c^{\pm}(\varphi_{\chi})$ by 
$u^{\pm}(\varphi_{\chi})$; see the paragraph after Theorem~\ref{thm:shimura}.)

The remaining cases, when at least one of $m,r$ or $k$ is odd, are absolutely 
similar. 
Likewise, the case of $L_f(m, {\rm Sym}^{2r} \varphi_{\chi})$, with its 
eight subcases depending on the parities of $m,r$ and $k$, is again very similar.
We leave the details to the reader.
\end{proof}

\section{Remarks on symmetric fourth power $L$-functions}
\label{sec:sym4}

\subsection{Consequences of cuspidality of the symmetric fourth}

We recall one of the main theorems of Kim--Shahidi \cite{kim-shahidi-duke} which 
characterizes cuspidality of the symmetric fourth power transfer of a cusp form on ${\rm 
GL}_2$. The 
following theorem is equivalent to \cite[Theorem 3.3.7]{kim-shahidi-duke}.

\begin{thm}
\label{thm:kim-shahidi-duke}
Let $F$ be a number field and let $\pi$ be a cuspidal automorphic representation of 
${\rm GL}_2({\mathbb A}_F)$. Then ${\rm Sym}^4(\pi)$ is not cuspidal as an automorphic 
representation of ${\rm GL}_5({\mathbb A}_F)$ if and only if 
$\pi$ is dihedral, tetrahedral or octahedral. 
\end{thm}

We digress a little and clarify the various equivalent versions
of a cuspidal automorphic representation being dihedral, tetrahedral or
octahedral type. The proofs are slightly scattered over the
literature and the aim is to guide the reader to the appropriate
references, while sketching some easy arguments. The expert on these issues
can skip to the paragraph after the proof of Proposition~\ref{prop:octahedral}.

\begin{prop}
\label{prop:dihedral}
Let $F$ be a number field and let $\pi$ be a cuspidal automorphic
representation of ${\rm GL}_2({\mathbb A}_F)$. 
Then the following are
equivalent.
\begin{enumerate}
\item $\pi$ is dihedral, i.e., $\pi$ is attached to a two dimensional 
irreducible representation $\sigma$ of $W_F$ which is induced from a 
character of $W_K$ for a quadratic extension $K/F$; 
$\sigma = {\rm Ind}_{W_K}^{W_F}(\chi)$. 
\item There exists a nontrivial id\`ele class character $\eta$ of $F$ 
such that $\pi \simeq \pi\otimes\eta$. The character $\eta$ is necessarily
quadratic. 
\item ${\rm Sym}^2(\pi)$ is not cuspidal. 
\end{enumerate}
\end{prop}

\begin{proof}
For (1) $\Rightarrow$ (2) take $\eta = \omega_{K/F}$. 
The statement $(2) \Rightarrow (1)$ is due to  
Labesse--Langlands \cite[Proposition 6.5]{labesse-langlands}.
The statement (2) $\Rightarrow$ (3)
can be seen using $L$-functions and the heuristic
$\sigma \otimes \sigma^{\vee} \simeq 
(\sigma \otimes \eta) \otimes \sigma^{\vee} \simeq
({\rm Sym}^2(\sigma)\otimes {\rm det}(\sigma)^{-1}\eta) \oplus \eta.$
Hence the $L$-function $L(s, {\rm Sym}^2(\pi) 
\otimes \omega_{\pi}^{-1}\eta)$
has a pole at $s=1$.
The statement (3) $\Rightarrow$ (2)  is contained in
Gelbart--Jacquet \cite[Theorem 9.3]{gelbart-jacquet}. 
\end{proof}

\begin{prop}
\label{prop:tetrahedral}
Let $F$ be a number field and let $\pi$ be a cuspidal automorphic 
representation of ${\rm GL}_2({\mathbb A}_F)$. Assume that $\pi$ is not
dihedral. Then the following are equivalent.  
\begin{enumerate}
\item $\pi$ is tetrahedral, i.e., $\pi$ is attached to a two dimensional
irreducible representation $\sigma: W_F \to {\rm GL}_2({\mathbb C})$ 
whose image in ${\rm PGL}_2({\mathbb C})$ is isomorphic to $A_4$. 
\item There exists a nontrivial id\`ele class character $\eta$ of $F$
such that ${\rm Sym}^2(\pi) \simeq {\rm Sym}^2(\pi)\otimes\eta$. The 
character $\eta$ is necessarily cubic.
\item ${\rm Sym}^2(\pi)$ is cuspidal and ${\rm Sym}^3(\pi)$ is not cuspidal.
\end{enumerate}
\end{prop}

\begin{proof}
The statement (1) $\Rightarrow$ (2) follows from the representation theory 
of $A_4$; it suffices to observe that there is a character $\eta$ such that
${\rm Sym}^2(\sigma) \simeq {\rm Sym}^2(\sigma)\otimes \eta$. 

The statement (2) $\Rightarrow$ (1) can be seen as follows: 
Let $E/F$ be the cyclic extension of degree $3$ defined by $\eta$. Since $\pi$ satisfies 
${\rm Sym}^2(\sigma) \simeq {\rm Sym}^2(\sigma)\otimes \eta$, we get using
the first paragraph of the proof of Lemma 9.2 in \cite{shahidi-crm} that
the base change $\pi_E$ of $\pi$ to $E$ is monomial, i.e., 
$\pi_E = \pi(\sigma_E)$ for a representation $\sigma_E$ of $W_E$ 
which is induced from a quadratic extension $K/E$ as 
$\sigma_E = {\rm Ind}_{W_K}^{W_E}(\chi)$ for some character $\chi$.
Since $\pi_E$ is ${\rm Gal}(E/F)$-invariant, so is $\sigma_E$. Hence 
$\sigma_E$ extends to a representation, say $\sigma$ of $W_F$. That it 
extends may be seen by Lemma 7.9 of \cite{knapp-stein}; or by appealing to 
the fact that $H^2(W_F, {\mathbb C}^*)$ is trivial. Consider the image 
$I$ of $\sigma$ in ${\rm PGL}_2({\mathbb C})$. It is not cyclic because 
$\pi$ is cuspidal; it is not dihedral because $\pi$ is assumed not to be 
dihedral; it is not $S_4$ or $A_5$ because neither has an 
index $3$ subgroup; hence the image $I$ is $A_4$ or that $\sigma$ is 
tetrahedral type. Now consider $\pi(\sigma)$ and its base change 
$\pi(\sigma)_E$ to $E$. We have 
$\pi(\sigma)_E = \pi(\sigma|_{W_E}) = \pi(\sigma_E) = \pi_E$. Hence 
for some $i = 0, 1,2$ we have $\pi = \pi(\sigma) \otimes \eta^i = 
\pi(\sigma\otimes \eta^i)$. 

The equivalence (2) $\Leftrightarrow $ (3) is contained in 
\cite[Proposition 6.3]{kim-shahidi-annals}
\end{proof}

\begin{prop}
\label{prop:octahedral}
Let $F$ be a number field and let $\pi$ be a cuspidal automorphic
representation of ${\rm GL}_2({\mathbb A}_F)$. Assume that $\pi$ is neither
dihedral nor tetrahedral. Then the following are equivalent.
\begin{enumerate}
\item $\pi$ is octahedral, i.e., $\pi$ is attached to a two dimensional
irreducible representation $\sigma: W_F \to {\rm GL}_2({\mathbb C})$
whose image in ${\rm PGL}_2({\mathbb C})$ is isomorphic to $S_4$.
\item There exists a quadratic extension $E/F$ and 
there exists a nontrivial id\`ele class character $\eta$ of $E$
such that ${\rm BC}_{E/F}({\rm Sym}^2(\pi)) \simeq 
{\rm BC}_{E/F}({\rm Sym}^2(\pi))\otimes\eta$. The
character $\eta$ is necessarily cubic.
\item ${\rm Sym}^2(\pi)$ and ${\rm Sym}^3(\pi)$ are cuspidal, and ${\rm Sym}^4(\pi)$ is not cuspidal.
\end{enumerate}
\end{prop}

\begin{proof}
The statement $(1) \Rightarrow (2)$ follows from 
Proposition~\ref{prop:tetrahedral} by noting that $A_4$ is a normal subgroup
of index $2$ in $S_4$ and that ${\rm BC}_{E/F}$ commutes with ${\rm Sym}^2$
which can be seen by verifying it locally everywhere. 
The statement $(2) \Rightarrow (1)$ is
proved in \cite[Proposition 3.3.8 (2)]{kim-shahidi-duke} and is similar to the 
the proof of $(2) \Rightarrow (1)$ of Proposition~\ref{prop:tetrahedral} above.
The equivalence $(2) \Leftrightarrow (3)$ is contained in 
\cite[Proposition 3.3.6]{kim-shahidi-duke}. 
\end{proof}

Now let us consider Deligne's conjecture for the special values of 
$L(s, {\rm Sym}^4 \varphi)$ where $\varphi$ is a primitive form in 
$S_k(N,\omega)$. The $L$-function is, up to shifting by $2(k-1)$, the
standard $L$-function $L(s, {\rm Sym}^4(\pi(\varphi))$. If the 
representation ${\rm Sym}^4(\pi(\varphi))$ is not cuspidal, then 
appealing to the above cuspidality theorem we know that 
$\pi(\varphi)$ is either dihedral, tetrahedral or octahedral. 
In the dihedral case we have given a proof in \S\ref{sec:dihedral}. 
The proof in the tetrahedral and octahedral cases, 
if one may use the word proof in such a context, 
boils down to showing that there are no critical integers for 
$L(s, {\rm Sym}^4, \pi)$ and hence Deligne's conjecture is vacuously true!
The following well known lemma says that in these cases the cusp form we begin with is 
necessarily of weight one ($k=1$), and so from 
Lemma~\ref{lem:critical-sym-odd} and Lemma~\ref{lem:critical-sym-even}
it follows that there are no critical points.

\begin{lemma}
\label{lem:polyhedral}
Let $\varphi$ be a primitive form in $S_k(N,\omega)$. Let $\pi = \pi(\varphi)$ be the
cuspidal representation of ${\rm GL}_2({\mathbb A}_{\mathbb Q})$ attached to 
$\varphi$. If $\pi$ corresponds
to a two dimensional representation $\sigma$ of $W_{\mathbb Q}$ whose image in 
${\rm PGL}_2({\mathbb C})$ is finite, then $k=1$. 
(In particular, if $\pi(\varphi)$ is tetrahedral or
octahedral, then the modular form $\varphi$ we begin with is necessarily a
weight one form.)
\end{lemma}

\begin{proof}
If $k \geq 2$, then the image of 
$\sigma_{\infty} = {\rm Ind}_{{\mathbb C}^*}^{W_{\mathbb R}}(\chi_{k-1})$
in ${\rm PGL}_2({\mathbb C})$ is infinite. 
\end{proof}

We need to consider the case now when ${\rm Sym}^4(\pi(\varphi))$ is cuspidal. 
A possible strategy then is to appeal to the work of Mahnkopf \cite{mahnkopf2}
and apply his results on the special values of standard $L$-functions of ${\rm GL}_n$ 
to the particular case of $L(s, {\rm Sym}^4(\pi(\varphi))$. Recall that an important
part of the hypothesis in his work is that representation one begins with is 
cohomological. Since this is of independent interest, we consider this in the 
next subsection.

\subsection{Cohomological criterion}

In this subsection we recall the following theorem,
essentially due to Labesse-Schwermer \cite{labesse-schwermer},
which says that symmetric power lifts of a 
holomorphic modular form are cohomological. 

\begin{thm}
\label{thm:symmetric-cohomology}
Let $\varphi \in S_k(N,\omega)$ with $k \geq 2$. 
Let $n \geq 1$. Assume that
${\rm Sym}^n(\pi(\varphi))$ is a cuspidal representation of
${\rm GL}_{n+1}({\mathbb A}_{\mathbb Q})$. Let 
$$
\Pi = {\rm Sym}^n(\pi(\varphi)) \otimes \xi \otimes |\!|\ |\!|^s
$$
where $\xi$ is any id\`ele class character such that 
$\xi_{\infty} = {\rm sgn}^{\epsilon}$, with $\epsilon \in \{0,1\}$, 
and $|\!|\ |\!|$ is the ad\`elic norm. We suppose that $s$ and $\epsilon$ 
satisfy: 
\begin{enumerate}
\item If $n$ is even, then let $s \in {\mathbb Z}$ and 
$\epsilon \equiv n(k-1)/2 \pmod{2}$. 
\item If $n$ is odd then, we let $s \in {\mathbb Z}$ if k is even,
and we let $s \in 1/2+{\mathbb Z}$ if k is odd. We impose no condition 
on $\epsilon$. 
\end{enumerate}
Then 
$\Pi \in {\rm Coh}(G_{n+1}, \mu^{\vee})$ 
where $\mu \in X^+(T_{n+1})$ is given by 
$$
\mu = \left(
\frac{n(k-2)}{2}+s, \frac{(n-2)(k-2)}{2}+s,\dots, \frac{-n(k-2)}{2}+s
\right) = (k-2)\rho_{n+1}+s.
$$
Here $\rho_{n+1}$ is half the sum of positive roots of ${\rm GL}_{n+1}$.
In other words, the representation 
${\rm Sym}^n(\pi(\varphi)) \otimes \xi \otimes |\cdot |^s$, with $\xi$ and 
$s$ as above, contributes to 
cuspidal cohomology of the locally symmetric space 
${\rm GL}_{n+1}({\mathbb Q})\backslash {\rm GL}_{n+1}({\mathbb 
A}_{\mathbb Q})/
K_f K_{n+1, \infty}^{\circ}$ with coefficients in the local 
system determined by 
$\rho_{\mu^{\vee}}$, where $K_f$ is a deep enough open compact subgroup of 
${\rm GL}_{n+1}({\mathbb A}_{{\mathbb Q},f})$. (Here ${\mathbb 
A}_{{\mathbb Q},f}$ denotes the finite ad\`eles of ${\mathbb Q}$.) 
\end{thm}

\begin{proof}
See \cite[Theorem 5.5]{raghuram-shahidi}. 
\end{proof}

\begin{cor}
Let $\varphi \in S_k(N,\omega)$. Assume that $k \geq 2$ and that $\varphi$ is not dihedral. 
Then, up to twisiting by a quadratic character, ${\rm Sym}^n(\pi(\varphi))$
for $n=2,3,4$, contributes to cuspidal cohomology. 
\end{cor}

\subsection{Special values}
\label{sec:sym4-specialvalues}

As in the hypothesis of the above corollary, consider a holomorphic primitive modular form
$\varphi \in S_k(N,\omega)$. Assume that $k \geq 2$ and that $\varphi$ is not dihedral (the other
cases being done with, as far as special values are concerned). Let 
$\Pi = {\rm Sym}^4(\pi(\varphi))$. Then by Theorem~\ref{thm:symmetric-cohomology} we have 
$\Pi \in {\rm Coh}({\rm GL}_5, \mu^{\vee})$ where $\mu = (k-2)\rho_5$. We may appeal to 
Mahnkopf \cite{mahnkopf2}
and get information on the special values of $L(s,\Pi)$ and hence about
$L(s, {\rm Sym}^4\varphi)$. The purpose of this section is to record, what according to us, 
are some impediments of this strategy. 

\begin{enumerate}
\item {\it Nonvanishing hypothesis.} Mahnkopf's work \cite{mahnkopf2} is based on a 
certain
nonvanishing hypothesis. This hypothesis shows up in several other works on 
special values of $L$-functions which are based on cohomological interpretations 
of zeta integrals. See \cite[\S 6.2]{raghuram-shahidi} for a summary of the main
results of \cite{mahnkopf2} and this nonvanishing hypothesis. Eliminating this 
hypothesis, a problem which concerns archimedean zeta integrals, is an important 
technical problem. 

\item {\it Auxiliary twisting.} In \cite{mahnkopf2} there is an auxiliary character 
$\eta$ which has been brought in to finesse the bad places for the 
representations at hand. The presence of this character only gives the special 
values of certain twisted $L$-functions, and not any particular $L$-function that 
one might care about.  We believe that it is possible to work through Mahnkopf 
\cite{mahnkopf2}, while using the observation that special values of local 
$L$-functions are always rational. (See \cite[Lemme 4.6]{clozel2}.) This is work in 
progress, and we hope to report on this on a future ocassion.

\item {\it Explicit comparison of periods.} 
This is a far more philosophical problem. The periods of Harder and Mahnkopf (and more
generally those in \S\ref{sec:periods}) come by comparing 
rational structures on the Whittaker model and on a certain cohomology space. 
However, in Deligne's conjectures, the periods come by comparing rational structures
on the de Rham and Betti realization of the underlying motive. There is no obvious
comparison between these periods. Untill this problem is explicitly
solved, the best one can hope is to prove a theorem which only formally looks like 
the predictions made by Deligne \cite{deligne}. 
\end{enumerate}

\bigskip

\noindent {\it Address:} 

\noindent 
A.~Raghuram, Department of Mathematics, 
Oklahoma State University, 
401 Mathematical Sciences, 
Stillwater, OK 74078, USA. 
E-mail : {\tt araghur@math.okstate.edu}

\smallskip

\noindent
Freydoon~Shahidi, Department of Mathematics, 
Purdue University, 150 N. University Street, 
West Lafayette, IN 47907, USA. 
E-mail : {\tt shahidi@math.purdue.edu}

\end{document}